\documentclass[12pt]{article}

\usepackage{amsmath,amsfonts,amsthm,epsfig,latexsym,graphicx,amssymb}
\usepackage{amssymb,amsthm}
\usepackage{psfrag}

\newtheorem{Que}{Question}
\newtheorem{Cnj}{Conjecture}
\newtheorem{Thm}{Theorem}
\newtheorem{Corl}{Corollary}
\newtheorem{Rema}{Remark}

\newtheorem{Lm}{Lemma}[section]
\newtheorem{Def}[Lm]{Definition}
\newtheorem{Prop}[Lm]{Proposition}
\newtheorem{Rem}[Lm]{Remark}
\newtheorem{Cor}[Lm]{Corollary}
\newtheorem{Cla}[Lm]{Claim}

\def\bdef{\begin{Def}}
\def\endef{\end{Def}}

\def\bthm{\begin{Thm}}
\def\bcnj{\begin{Cnj}}
\def\ethm{\end{Thm}}
\def\ecnj{\end{Cnj}}
\def\bque{\begin{Que}}
\def\eque{\end{Que}}

\def\brema{\begin{Rema}}
\def\erema{\end{Rema}}

\def\bpro{\begin{Prop}}
\def\epro{\end{Prop}}

\def\blm{\begin{Lm}}
\def\elm{\end{Lm}}

\def\bcorl{\begin{Corl}}
\def\ecorl{\end{Corl}}

\def\bcor{\begin{Cor}}
\def\ecor{\end{Cor}}

\def\brm{\begin{Rem}}
\def\erm{\end{Rem}}

\def\bcl{\begin{Cla}}
\def\ecl{\end{Cla}}

\def\bfig{\begin{picture}}
\def\efig{\end{picture}}

\def\beq{\begin{eqnarray}}
\def\eneq{\end{eqnarray}}

\def\beal{\begin{aligned}}
\def\enal{\end{aligned}}

\def\neig{{V_0}}

 \def\NN{{\mathbb N}} 

 \def\RR{{\mathbb R}}

\def\mfT{{\mathfrak{T}}}

\def\Si{\Sigma}
\def\La{\Lambda}
\def\De{\Delta}
\def\Om{\Omega}

\def\la{\lambda}
\def\be{\beta}

\def\ga{\gamma}
\def\const{\Theta}

\def\ve{\varepsilon}

\def\diffM{\mbox{{\rm Diff\,}}^1(M)}

  \def\cG{{\cal G}}  \def\cS{{\cal
S}}
  \def\cH{{\cal H}}  
\def\cC{{\cal C}}  \def\cI{{\cal I}} \def\cO{{\cal O}} \def\cU{{\cal
U}}
  \def\cJ{{\cal J}}  \def\cV{{\cal
V}}
\def\cE{{\cal E}}    
   \def\cR{{\cal R}} 
\def\cY{{\cal Y}}

\def\diff{\operatorname{Diff\,}^1}

\def\rj{{\operatorname{j}}}
\def\ri{{\operatorname{i}}}
\def\st{{\operatorname{s}}}
\def\sst{\operatorname{ss}}
\def\cst{\operatorname{cs}}
\def\cent{\operatorname{c}}
\def\ct{{\operatorname{c}}}
\def\ut{{\operatorname{u}}}
\def\uut{\operatorname{uu}}

\def\loc{{\operatorname{loc}}}

\title{Non-hyperbolic ergodic measures for non-hyperbolic
homoclinic classes}
\author{Lorenzo J. D\'\i az and Anton Gorodetski}

\makeatletter

\def\keywords#1{{\def\@thefnmark{\relax}\@footnotetext{#1}}}

\let\subjclass\keywords

\makeatother

\begin{document}

\maketitle






\begin{abstract}
We prove that there is a residual subset $\mathcal{S}$ in
$\text{Diff\,}^1(M)$ such that, for every $f\in \mathcal{S}$, any
homoclinic class of $f$ containing saddles of different indices
(dimension of the unstable bundle) contains also an uncountable
support of an invariant ergodic non-hyperbolic (one of the Lyapunov
exponents is equal to zero) measure of $f$.
\end{abstract}


\section{Introduction}

 It was shown in the 1960s by Abraham and  Smale
that uniform hyperbolicity is not dense in the space of dynamical
systems \cite{AS}. This necessitated weakening the notion of
hyperbolicity. One of the possible approaches (Pesin's theory
\cite{Pe}) is to characterize hyperbolic behavior by non-zero
Lyapunov exponents with respect to some invariant measure. The most
natural case is that of a system with a smooth invariant measure. In
this setting Lyapunov exponents were studied in  various aspects,
such that removability of zero exponents \cite{BB,SW}, genericity of
zero or non-zero exponents \cite{B,BcV}, and existence of hyperbolic
measures \cite{DP}.   However, the question about Lyapunov exponents
can also be considered for non-conservative maps. 

Recall that if $\mu$ is an ergodic measure of a diffeomorphism $f:
M\to M,\, \text{\rm dim}\, M={\mathbf{m}},$ then there is a set
$\Lambda$ of full $\mu$-measure and real numbers  $\chi^1_{\mu}\le
\chi^2_{\mu}\le \dots \le \chi^{\mathbf{m}}_{\mu}$ such that, for
every $x\in \Lambda$ and every non-zero vector $v\in T_xM$, one
has $\lim_{n\to \infty} (1/n) \, \log || Df^n(v)|| = \chi^i_\mu$
for some $i=1, \dots, \mathbf{m}$, see \cite{O, M2}. The number
$\chi^i_\mu$ is  the \emph{$i$-th Lyapunov exponent} of the
measure $\mu$. \bdef An ergodic invariant measure of a
diffeomorphism is called {\it non-hyperbolic} if at least one of
its Lyapunov exponents is equal to zero.
\endef

Some natural questions arise while considering non-hyperbolic
diffeomorphisms and Lyapunov exponents of their ergodic measures.
First, does a generic diffeomorphism have non-zero Lyapunov exponents
for each invariant measure?  A negative answer to this question was
given recently in \cite{KN} by constructing a $C^1$-open subset in
$\text{Diff\,}^r(M)$ ($r\ge 1, \,$ $M$ is a closed smooth manifold,
$ \text{dim}\, M\ge 3$) of diffeomorphisms exhibiting non-hyperbolic
ergodic invariant measures with uncountable support.

It seems interesting  also to consider non-hyperbolic invariant
measures with respect to other questions. How to characterize the
absence of uniform hyperbolicity? What dynamical structures can not
exist in the uniformly hyperbolic setting but must be present in the
complement? A number of conjectures related to this question had
been stated. The most influential one is the Palis' Conjecture 
\cite{Pa} that
claims, roughly speaking, that diffeomorphisms exhibiting homoclinic
tangencies or heterodimensional cycles are dense in the complement
to the set of uniformly hyperbolic systems. This conjecture was
proved  in \cite{PS} for $C^1$ surface diffeomorphisms (in that
case, heterodimensional cycles are not considered). Another
candidates for a ``non-hyperbolic structure" are, for instance,
super-exponential growth of the number of periodic points
\cite{K,BDF} and the absence of shadowing property
\cite{BDT,YY,AD,S}. Here we would like to suggest a reformulation of
Palis' conjecture meaning that the non-hyperbolic behavior is
detected in the ergodic level:

\bcnj \label{c.dichotomy}  In $\text{{\rm Diff\,}}^r(M), r\ge 1,$
there exists an open and dense subset $\mathcal{U}\subset
\text{{\rm Diff\,}}^r(M)$ such that every diffeomorphism $f\in
\mathcal{U}$ is either uniformly hyperbolic or has an ergodic
non-hyperbolic invariant measure. \ecnj

We observe that this conjecture holds if we replace the \emph{open
and dense} condition by just \emph{dense}.\footnote{A
$C^1$-diffeomorphism $f$ satisfies the star condition if it has a
$C^1$-neighborhood $\mathcal{U}$ such that every periodic point of
every diffeomorphism in $\mathcal{U}$ is hyperbolic. The star
condition is equivalent to the Axiom A and no-cycles conditions,
see \cite{A,H,L,M1}. Thus to prove the weak ``conjecture" it
suffices to approximate a non-star diffeomorphism by a
diffeomorphism with a non-hyperbolic periodic point and consider a
measure supported on that orbit.} Alternative (weaker) slightly
different reformulation is to consider generic diffeomorphisms.

We observe that, as a consequence of our main result, this
conjecture holds in the so-called tame setting, see
Theorem~\ref{t.main} and the discussion below.

\brm A shift along the phase curves of the Bowen's example of an
invariant disk bounded by the separatrices of two  saddles provides
an example of a non-Axiom A diffeomorphism without non-hyperbolic
invariant measures, see {\em \cite{T}}. For another
less degenerate example, see {\em\cite{CLR}}. \erm

\brm If all Lyapunov exponents for every invariant measure of a
local $C^1$-diffeomorphism  are positive then it is uniformly
expanding. Also, if a partially hyperbolic diffeomorphism has all
central Lyapunov exponents positive for every invariant measure then
it is uniformly hyperbolic. Note that this hypothesis necessarily
implies that all periodic points are hyperbolic and have the same
index. See {\em\cite{AAS}} for both results. In this paper, we study
a different setting: we consider transitive sets (homoclinic
classes)  containing saddles of different indices and construct
non-hyperbolic ergodic measures. \erm

Let us recall that the \emph{homoclinic class} of a hyperbolic
periodic  point $P$ of a diffeomorphism $f$, denoted by $H(P, f)$,
is the closure of the transverse intersections of the invariant
manifolds (stable and unstable ones) of the orbit of $P$ (note
that the homoclinic class of a sink or a source is just its
orbit). A homoclinic class can be also defined as the closure of
the set of hyperbolic saddles $Q$ \emph{homoclinically related} to
$P$ (the stable manifold of the orbit of $Q$ transversely meets
the unstable one of the orbit of $P$ and vice-versa). Note that
two homoclinically related saddles have the same \emph{index}
(dimension of the unstable bundle).
 In \cite{N2}  the notion
of a homoclinic class was proposed as a generalization of a
uniformly hyperbolic basic set of an Axiom A diffeomorphism.
However, homoclinic classes may fail, in general, to be hyperbolic
and may contain saddles having different indices (in fact, this is
the context of our paper). For explicit nontrivial examples of
non-hyperbolic homoclinic classes see \cite{BD1,D} (see also
\cite{GI1,G} for similar examples studied later by different
methods).

In order to support Conjecture \ref{c.dichotomy}, we prove the
following theorem.

\bthm\label{t.main} In $\text{\rm{Diff\,}}^1(M)$ there exists a
residual subset $\mathcal{S}$ such that, for every $f\in
\mathcal{S}$, any homoclinic class containing saddles of different
indices contains also an uncountable support of an invariant
ergodic non-hyperbolic measure of $f$. \ethm

In fact, we prove more:

\
\newline
{\bf Addendum.} {\it Given a diffeomorphism $f\in \cS$, every
homoclinic class of $f$ with saddles whose stable bundles have
dimensions $s$ and $s+r$, $r\ge 1$, for every $j=s+1,\dots,s+r$
contains  a uncountable support of an invariant ergodic measure
whose $j$-th Lyapunov exponent is
zero.} 

\medskip

We also would like to mention that recently the following result
had been announced in \cite{ABC}: for $C^1$-generic
diffeomorphisms\footnote{By $C^1$-generic diffeomorphisms we mean
diffeomorphisms in a residual subset of $\text{Diff\,}^1(M)$.} the
generic measures supported on isolated homoclinic classes are
ergodic and hyperbolic (all Lyapunov exponents are non-zero).

\medskip

Let us  now give some motivation of our result. First we quote the question posed by Shub and Wilkinson.
\bque\label{q.sw} {\rm (\cite[Question 2]{SW})} For $r \ge 1$, is it true for the generic
$f\in\text{{\rm Diff\,}}^r(M)$ and any weak limit $\nu$ of
averages of the push forwards
$\frac{1}{n}\sum_{i=1}^{n}f_*^i\text{\it Leb}$ that almost every
ergodic component of $\nu$ has some exponents not equal to 0
($\nu$-a.e.)? All exponents not equal to 0? \eque

Even if one considers all the invariant measures, not only limit
points of the averages of shifts of the Lebesgue measure,  the
question remains nontrivial and meaningful.\footnote{Numerically presence of zero Lyapunov exponents was studied in \cite{GOST}, where some zero Lyapunov exponents were obtained.  Whether this numerical effect is really
related to the presence of non-hyperbolic measures, or it is an
artifact of numerical computations, this is not clear so far.}

Theorem~\ref{t.main}  was also motivated by the construction in
\cite{GIKN} of non-hy\-per\-bo\-lic ergodic measures for a class
of
 skew products $F\colon \Sigma\times
\mathbb{S}^1\to \Sigma\times \mathbb{S}^1$ of the form
$F(\omega,x)=(\sigma(\omega), f_{\omega_0}(x))$, where
$\Sigma=\{0,1\}^{\mathbb{Z}}$, $\sigma$ is the shift in $\Sigma$,
and $f_0$ and $f_1$ are appropriate circle diffeomorphisms. We use the method developed in \cite{GIKN} for
obtaining non-hyperbolic ergodic measures as weak limits of measures
supported on periodic points.

It is also related to the series of results in
\cite{ABCDW,BC,BDF,BDP,CMP} about the geometrical structure of
non-hyperbolic homoclinic classes of $C^1$-generic
diffeomorphisms.
 In a sense our
results give a description of the dynamics of  non-hyperbolic
homoclinic classes in the ergodic level. A key fact here is that
these classes exhibit \emph{heterodimensional cycles} in a
persistent way. The analysis of the dynamics of these cycles is
another ingredient in this paper. Recall that a diffeomorphism $f$
has a \emph{  heterodimensional cycle} if there are saddles $P$
and $Q$ of $f$ having different indices such that their invariant
manifolds meet cyclically (i.e., $W^s(P,f)\cap
W^u(Q,f)\ne\emptyset$ and $W^u(P,f)\cap W^s(Q,f)\ne\emptyset$).

\medskip

Let us suggest a naive way to prove Theorem~\ref{t.main}. By the
results in \cite{ABCDW}, there is a sequence of saddles in the
homoclinic class with a central Lyapunov exponent converging to
zero. One is tempted to consider a weak limit of the invariant
atomic measures supported on those orbits and to expect that the
resulting measure is non-hyperbolic. Unfortunately,  the resulting
measure, for instance, can be supported on several hyperbolic
periodic orbits. To get a non-hyperbolic ergodic measure, we need
to choose those periodic orbits in an intricate  way to guarantee
that the limit measure is ergodic. We use here the strategy
suggested by Ilyashenko, see \cite{GIKN}. Roughly speaking, we
construct a sequence of periodic orbits such that each of the
orbits shadows the previous one for a long time, but differs from
it for a much shorter time. This provides simultaneously
decreasing Lyapunov exponents and ergodicity of the limit measure.
To generate such a sequence of orbits we use heterodimensional
cycles.

\medskip

We observe that the unfolding of any co-index one cycle (a cycle
has {\emph{co-index one}} if $\mbox{index}(P)=\mbox{index }(Q)\pm
1$), say associated to $f$, generates an open set of
$C^1$-diffeomorphisms $\cO$ such that $f$ is in the closure of
$\cO$ and  every $g$ in a residual subset $\cS$ of $\cO$ has two
saddles $A_h$ and $B_h$ having different indices such that
$H(A_h,h)=H(B_h,h)$. This result is a consequence of the
constructions in \cite{BD4}\footnote{This result is not  stated
explicitly in \cite{BD4}, where the results are stated in terms of
$C^1$-robust cycles, but it follows immediately from the
``blender-like" constructions there.}. Thus we can apply
Theorem~\ref{t.main} to the set $\cO$ getting the following:

\bcorl Let $f$ be a $C^1$-diffeomorphism with a co-index one cycle. Then
there are a $C^1$-open set $\,\cO$ and a residual subset $\,\cR$ of
$\,\cO$ such that $f$ is in the closure of $\,\cO$ and every $g\in
\cR$ has a homoclinic class containing the uncountable support of
a
 non-hyperbolic ergodic measure.
\ecorl


\subsubsection*{Open questions and consequences.}

First of all,  we observe that in Theorem~\ref{t.main} the support
of the non-hyperbolic measure is in general properly contained in
the homoclinic class. To construct this measure a key ingredient is
\emph{partial hyperbolicity} (with one dimensional central
direction),  and we need to identify a part
of the homoclinic class where the \emph{relative
dynamics} is partially hyperbolic. We consider measures having
supports contained in that partially hyperbolic region.

On the one hand, these comments imply that, for homoclinic classes
whose non-hyperbolic central bundle has dimension equal to or
greater than two,
 our method do not provide non-hyperbolic measures with full support in
the class (see, for instance, \cite{BnV} for examples of diffeomorphisms
with a homoclinic class equal to the whole ambient manifold whose
non-hyperbolic central bundle has dimension two).
 On the other hand, for homoclinic classes
with a partially hyperbolic splitting with one-dimensional central
bundle, one can expect to extend our method to construct
non-hyperbolic ergodic measures with full support. However, a
first  difficulty for such an extension is to understand the
distribution of the periodic orbits in the class. A first step in
this direction is the result in \cite{C} claiming that
$C^1$-generic homoclinic classes  are Hausdorff limits of periodic
orbits. These comments lead to the following question:

\bque \label{q.fullsupport} Consider a homoclinic class $H(P,f)$ of a $C^1$-generic diffeomorphism $f$
containing saddles of different indices and having a partially hyperbolic splitting with one-dimensional central
direction. Is there a non-hyperbolic ergodic measure whose support is the whole class $H(P,f)$\footnote{In the
recent work \cite{Na} it was shown, in particular, that in the construction in \cite{KN} one can find a
non-hyperbolic invariant measure supported on the whole homoclinic class.}?.
 \eque

A diffeomorphism $f$ is \emph{tame} if its  homoclinic classes are   robustly
isolated. Tame diffeomorphisms form an open set in
$\text{\rm{Diff\,}}^1(M)$.
 There is a residual subset
$\mathcal{R}$ of $\text{Diff\,}^1(M)$ such that for every $f\in
\mathcal{R}$ homoclinic classes of the saddles of $f$ form a
partition of a part of the limit set of $f$, see \cite{CMP}. The
set $\mathcal{R}$ is the union of two disjoint sets $\mathcal{T}$
and $\mathcal{W}$ which are relatively open in $\mathcal{R}$. The
set $\mathcal{T}$ (the intersection of the set of tame
diffeomorphisms and the residual set $\mathcal{R}$) coincides with
the set of diffeomorphisms with finitely many homoclinic classes
and $\mathcal{W}$ (the so called \emph{wild diffeomorphisms}) is
the set of diffeomorphisms with infinitely many homoclinic
classes. If the dimension of $M$ is strictly greater than $2$
both sets $\mathcal{T}$ and $\mathcal{W}$ are non-empty, see
\cite{BD2,BD3}.

For  tame diffeomorphism every homoclinic class is either hyperbolic
or contains saddles of different indices \cite{ABCDW}. In particular,  a
generic version of Conjecture~\ref{c.dichotomy} holds for tame
diffeomorphisms:

\bcorl\label{c.tame} Every generic tame diffeomorphism  is either
hyperbolic or has a non-hyperbolic invariant ergodic measure with
uncountable support.
 \ecorl
This paper shows that to settle a generic version of Conjecture~\ref{c.dichotomy} it is
enough to prove that every wild diffeomorphism has some homoclinic
class containing saddles having different indices. At present time,
this fact is not known, although for all known examples of wild
diffeomorphisms  this occurs. 

We also expect that the following natural generalization of Theorem
1 holds.

\bcnj \label{c.homoclinic} In $\text{{\rm Diff\,}}^1(M)$ there
exits a residual subset $\mathcal{R}\subset \text{{\rm
Diff\,}}^1(M)$ such that for every diffeomorphism $f\in
\mathcal{R}$ every homoclinic class is either uniformly hyperbolic
or contains a support of an ergodic non-hyperbolic invariant
measure. \ecnj

We observe that
recently the proof of this conjecture for generic diffeomorphisms far from homoclinic tangencies (using the results in our paper) was announced in \cite{Y}.

A diffeomorphism is {\emph{transitive}} if it has a dense orbit, and
is {\emph{$C^r$-robustly transitive}} if it has a $C^r$-neighborhood
consisting of transitive diffeomorphisms.
Most of the examples 
 mentioned above are partially hyperbolic transitive systems. In fact, $C^1$-robust transitivity implies some form of weak
hyperbolicity \cite{DPU,BDP}. This leads to the following weaker
version of Conjecture~\ref{c.dichotomy}.

\bcnj \label{c.robust} Denote by $\mathcal{RT}^{\,r}\subset
\text{\rm Diff\,}^r(M)$ the (open) set of robustly transitive
diffeomorphisms.  In $\mathcal{RT}^{\,r}$ there exists an open and
dense subset $\mathcal{U}\subset \mathcal{RT}^{\,r}$ such that
every diffeomorphism $f\in \mathcal{U}$ is either uniformly
hyperbolic or has an ergodic non-hyperbolic invariant measure.
\ecnj

Since open and densely  robustly transitive $C^1$-diffeomorphisms are
tame diffeomorphisms (they have just one homoclinic class equal to
the ambient manifold), Corollary~\ref{c.tame} implies that this
conjecture holds $C^1$-generically.

\medskip

In order to have homoclinic classes with saddles of different
indices we need to consider a phase space of dimension at least
three. Let us shortly comment the two-dimensional case.  In that
case  it is expected that Axiom A diffeomorphisms are dense in
$\text{Diff\,}^1(M)$ (although no proof have been given
yet)\footnote{For recent results about the $C^1$-density of Axiom
A surface diffeomorphisms  and a discussion of the current state
of this problem, see \cite{ABCD}.}. In $\text{Diff\,}^r(M),\ r\ge
2,$ there are {\emph{Newhouse domains}}, i.e., open sets where
there are dense subsets of diffeomorphisms exhibiting homoclinic
tangencies \cite{N1} and non-hyperbolic periodic points
\cite{GST}, and where generic $C^2$-diffeomorphisms have infinite
number of sinks or sources \cite{N3}.\footnote{In higher
dimensions there are $C^1$-open sets where the diffeomorphisms
with infinitely many sinks or sources are generic \cite{BD2}.}

\bque \label{e.newhouse} Let $\mathcal{U}\subset \text{\rm
Diff\,}^r(M),\ \text{dim}\ M=2, \ r\ge 2,$ be a Newhouse domain. Is
it true that a generic diffeomorphism from $\mathcal{U}$ has an
invariant ergodic non-hyperbolic measure? A similar question can be
posed in the $C^1$-topology for locally generic diffeomorphisms with
infinitely many sinks or sources. \eque

We observe that \cite{CLR} provides examples of homoclinic classes
of surface diffeomorphisms containing non-transverse intersections (tangencies) whose ergodic
measures have non-zero Lyapunov exponents. A positive answer to this
question should mean that such a situation is quite pathological.

\subsubsection*{Organization of the paper.}

This paper is organized as follows. In
Section~\ref{s.ergodicmeasures}, we review the construction in
\cite{GIKN} of
non-hyperbolic ergodic measures (with uncountable support)
obtained as limit of measures supported on periodic orbits.
Notice that a key ingredient of this construction is partial
hyperbolicity (with one dimensional central direction) in a
fixed region of the manifold.

Section~\ref{s.generic} consists of two parts.
Section~\ref{ss.generic} contains results about homoclinic classes
of $C^1$-generic diffeomorphisms. Using these results, we fix the
residual subset of $\diffM$ that we will consider in our
constructions. In Section~\ref{ss.creationofcycles}, we state a
perturbation result about generation of heterodimensional cycles
in homoclinic classes containing saddles of different indices
(non-hyperbolic classes).

In Section~\ref{s.generationzero}, for non-hyperbolic homoclinic
classes, we state Propositions~\ref{p.simpleanddense} and
\ref{p.bdf} about  generation of saddles with central Lyapunov
exponents close to zero and of persistent heterodimensional cycles
(involving such saddles). If a $C^1$-generic homoclinic class
contains saddles having different indices then we can generate a
co-index one cycle whose relative dynamics in a neighborhood $V$ of
the cycle is partially hyperbolic (with one dimensional central
direction). The unfoldings of these cycles generate saddles (whose
orbits are contained in the fixed neighborhood $V$) which are in the
homoclinic class that we consider and have central Lyapunov
exponents close to zero.  A key fact is that we can use such saddles
to get new partially hyperbolic heterodimensional cycles (in the set
$V$) and new saddles. In this way, we get sequences of saddles whose
central Lyapunov exponents go to zero.  We will apply to these
saddles the results in Section~\ref{s.ergodicmeasures} to get the
non-hyperbolic ergodic measures in Theorem~ \ref{t.main}. This is
done in Section~\ref{s.genzero} by using an inductive argument
combining the propositions above
 (generation of saddles and cycles)
and the results in Section~\ref{s.ergodicmeasures}.

To prove Propositions~\ref{p.simpleanddense} and \ref{p.bdf}  we need to adapt previous constructions
about cycles and homoclinic classes in \cite{ABCDW,BD4,BDF} to our setting. A difficulty here is that
we need to consider the relative dynamics in the fixed region $V$ above. We note that in the constructions
in previous papers the perturbations may be global ones and then much more general perturbation results can be applied
directly.

\subsubsection*{Standing notation}
\begin{itemize}
\item
By a perturbation of $f$ we always mean a diffeomorphism which is $C^1$ arbitrarily close to $f$.
\item
Given a hyperbolic periodic point $P_f$ of a diffeomorphism $f$,
for a $C^1$-close map $g$  we denote by $P_g$ the continuation of
the point $P_f$ for $g$.
\item
We will denote by $\pi(P_f)$ the period of a periodic point $P_f$ of $f$.
\item
Since there is no possibility of confusion, by  a {\emph{cycle}} we
always mean a heterodimensional cycle.
\item
 Partially hyperbolic sets that we consider are always
{\emph{strongly
 partially hyperbolic,}} that is, they have partially hyperbolic splittings
 with three subbundles, $T_xM=E_x^{\sst}\oplus E_x^{\cent}\oplus
 E^{\uut}_x$,
where $E^{\sst}$ is uniformly contracting, $E^{\uut}$ is uniformly
expanding, and the central subbundle $E^{\cent}$ is one
dimensional.
\end{itemize}

\subsubsection*{Acknowledgments.} {LJD
is partially supported by CNPq, Faperj, and PRONEX (Brazil). LJD thanks the kind hospitality and the
support of CALTECH and UC Irvine during his visits while preparing this paper. AG was partially supported by the grant RFFI-CNRS 050102801. The authors also thank the hospitality and the support
of the International Workshop on Global Dynamics Beyond Uniform Hyperbolicity (Northwestern University), and
International Symposium of Dynamical Systems (Federal University of Bahia). Also, we thank F.~Abdenur,
Ch.~Bonatti, Yu.~Ilyashenko, V.~Kleptsyn, M.~Nalsky, and A.~Wilkinson for useful discussions.}


\section{Periodic points and non-hyperbolic ergodic measures}
\label{s.ergodicmeasures}

In this part of the paper we follow the approach suggested in
\cite{GIKN} to provide sufficient conditions for the existence of
non-hyperbolic ergodic invariant measures with uncountable
support, see Proposition~\ref{p.firstmainproposition} below.


\subsection{Ergodicity, invariant direction fields, and Lyapunov exponents}
\label{ss.direction} Assume that a diffeomorphism $f:M\to M$ has
an invariant continuous direction field $E$ in an open set
$O\subset M$. Then for every invariant measure $\mu$ whose support
(denoted by $\rm{supp}\ \mu$) is contained in $O$ one of the
Lyapunov exponents of $\mu$ is associated to $E$ (denote it by
$\chi^E$). Namely, for $\mu$-almost every point $x\in M$ and for a
non-zero vector $v\in T_xM$ from the corresponding direction $E$,
$$
\lim_{n\to \infty}\frac{1}{n}\log|Df^n(v)|=\chi^E(\mu).
$$

Speaking about convergence of measures we always mean $*$-weak
convergence: $\mu_n$ converges to  $\mu$, if for any continuous
function $\varphi\colon M\to \mathbb{R}$ it holds
$$
\int \varphi \, d\mu_n \to \int \varphi \, d\mu, \qquad \mbox{as }
n\to\infty.
$$

\blm \label{l.limit} Let diffeomorphism $f:M\to M$ has an
invariant continuous direction field $E$ in an open set $O\subset
M$. Let $\mu_n$ and $\mu$ be ergodic probability measures with
supports in $O$, and $\mu_n\to \mu$ as $n\to \infty$. Then
$\chi^E(\mu_n)\to \chi^E(\mu)$. \elm

\begin{proof}
Define the continuous function $\varphi\colon O\to \mathbb{R},\,
\varphi(x)=\log \frac{|Df_x(v_x)|}{|v_x|}$, where $v_x\in T_x M$
is any non-zero vector from the direction $E$ depending
continuously on $x$. By definition, Lyapunov exponent $\chi^E$
along the direction field $E$ at $x$ is a time average of the
function $\varphi$ at this point. Due to ergodicity of measures
$\mu_n$ (respectively, $\mu$), this time average is equal to the
space average of the function $\varphi$ with respect to the
corresponding measure for $\mu_n$- (respectively, $\mu$-) almost
every point. Since the function $\varphi$ is continuous, due to
the $*$-weak convergence of measures $\mu_n \to \mu$, we have:
$$
\chi^E(\mu_n) = \int \varphi \, d\mu_n \to \int \varphi \, d\mu =
    \chi^E(\mu).
$$
This completes the proof of the lemma.
\end{proof}

\subsection{Sufficient conditions for ergodicity}
\label{ss.sufficientergodicity}

Let $\mathfrak{G}$ be an arbitrary continuous map of a metric
compact space $\mathfrak{Q}$ into itself. Assume that $X_n$ are
periodic orbits of the map $\mathfrak{G}$, $\pi(X_n)$ are their
periods, and $\mu_n$ are atomic measures uniformly distributed on
these orbits.

\bdef Let us call {\em $n$-measure of the point $x_0$} an atomic
measure uniformly distributed on $n$ subsequent iterations of the
point $x_0$ under the map $\mathfrak{G}$:
$$
\nu_n(x_0) = \frac{1}{n} \sum\limits_{i=0}^{n-1}
\delta_{\mathfrak{G}^i(x_0)},
$$
where $\delta_x$ is  $\delta$-measure supported at point $x$.
\endef

The next lemma is a key point in the proof of the ergodicity of a
limit measure.  This lemma 
was suggested by Yu.~Ilyashenko, who was inspired by the ideas of
the work by A.~Katok and A.~Stepin on periodic approximations of
ergodic systems, see~\cite{KS1,KS2}.

\blm\label{ergodicitylemma} {\bf (\cite[Lemma 2]{GIKN})} Let
$\{X_n\}$ be a sequence of periodic orbits with increasing periods
$\pi(X_n)$ of a continuous map $\mathfrak{G}$ of a metric compact
space  $\mathfrak{Q}$ into itself. For each $n$, let $\mu_n$ be
the probability atomic measure uniformly distributed on the orbit
$X_n$. 

Assume that for each continuous function $\varphi$ on
$\mathfrak{Q}$ and all
 $\varepsilon > 0$ there exists $N = N(\varepsilon, \varphi) \in \Bbb{N}$ such that  for all $m>N$
 there exists a subset $\widetilde{X}_{m,\varepsilon} \subset X_m$, which satisfies the following conditions:
\begin{enumerate}
\item $ \mu_m(\widetilde{X}_{m,\varepsilon}) > 1 - \varepsilon$
and
\item for any  $n$, such that $m>n\ge N$, and for all $x \in \widetilde{X}_{m,\varepsilon}$
$$
\left| \int\varphi \, d\nu_{\pi(X_n)}(x) - \int\varphi \, d\mu_n
\right| < \varepsilon.
$$
\end{enumerate}
Then every limit point $\mu$  of the sequence $\{\mu_n\}$ is an
ergodic measure. \elm

\subsection{Sufficient conditions for existence of an invariant non-hyperbolic measure}
\label{ss.sufficientnonhyperbolic}

Given a finite set $\Gamma$ denote by $\#\Gamma$ the cardinality
of $\Gamma$.

\bdef \label{d.goodaprox} A periodic orbit $Y$ of a map $f$ is a
$(\gamma, \varkappa)$-good approximation of the periodic orbit $X$
of $f$ if the following holds.
\begin{itemize}
  \item There exists a subset $\Gamma$ of \,$Y$
and a projection $\rho \colon \Gamma \to X$ such that
$$
\text{\rm dist} (f^j(y),f^j(\rho(y)))<\gamma,
$$
for every $y\in \Gamma$ and every $j=0,1,\dots, \pi(X)-1$;
 \item $\dfrac{\#\Gamma}{\#Y}\ge \varkappa$; 
  \item $\#\rho^{-1}(x)$ is the same for all $x\in X$.
\end{itemize}
\endef

\bpro \label{p.firstmainproposition} Assume that a diffeomorphism
$f:M\to M$ has the following properties:
\begin{description}
\item{{\bf 1)}}
there exists an open domain $O\subset M$ such that $f$ has
  an invariant continuous direction field $E$ in $O$;
\item{{\bf 2)}}
there exists a sequence of periodic orbits
  $\{X_n\}_{n=1}^{\infty}$ of $f$
whose periods $\pi(X_n)$
  tend to infinity as $n\to \infty$
and such that $
  \overline{\cup_{n=1}^{\infty}X_n}\subset
  O.$
\end{description}
   Denote by $\chi^E(X)$ the Lyapunov exponent of $f$ along the orbit $X$ with
respect to the invariant direction field $E$.
\begin{description}
\item{{\bf 3)}}
There exists a sequence of numbers $\{\gamma_{n}\}_{n=1}^{\infty},
\gamma_n>0$, and a constant $C>0$ such that for each $n$ the orbit
$X_{n+1}$ is a $(\gamma_n, 1-C\,|\chi^E(X_n)|)$-good approximation
of the orbit $X_{n}$;
\item{{\bf 4)}}
let $d_n$ be the minimal distance between the points of the orbit
$X_n$, then
$$
\displaystyle{\gamma_n< \frac{\min_{1\le i\le n} d_i}{3\cdot
2^n}};
$$
\item{{\bf 5)}} there exits a constant $\xi\in(0,1)$ such that for every $n$
$$|\chi^E(X_{n+1})|< \xi\, |\chi^E(X_n)|.$$
\end{description}
Then $f$ has a non-hyperbolic invariant ergodic measure with an
uncountable support.
  \epro

Of course, in the previous proposition, the obtained
non-hyperbolic invariant measure $\mu$ has zero Lyapunov exponent
along the direction $E$.  The measure $\mu$ is a weak limit point
of the sequence of measures supported in the orbits $X_n$.

\brm
In the context of this paper the claim that the support of the constructed non-hyperbolic measure is uncountable is not essential. Indeed, we consider $C^1$-generic diffeomorphisms whose periodic points are hyperbolic. But we provide this claim keeping
in mind future applications (considering open sets of diffeomorphisms).
\erm

  \begin{proof}
  Let $\mu_n$ be the probability
  atomic measure uniformly distributed on the orbit $X_n$. Let us
  check that the conditions of  Lemma \ref{ergodicitylemma} are
  satisfied by these measures.

  Set $\varkappa_n=1-C|\chi^E(X_n)|$.
Take arbitrary $\varepsilon > 0$ and continuous map $\varphi
\colon M \to \mathbb{R}$. By assumption 3), for orbits $\{X_n\}$ a
sequence of subsets $\widetilde{X}_n\subset X_n$ and a sequence of
projections $\rho_n: {\widetilde X}_{n+1} \to X_n$ are defined
such that:

\begin{equation}
\label{e.varkappa} \prod\limits_{n=1}^{\infty} \dfrac{\#
{\widetilde X}_{n+1} } {\# X_{n+1}}\ge \prod\limits_{n=1}^{\infty}
\varkappa_n =\hat \varkappa
     > 0.
     \end{equation}
Indeed, the product is convergent (and different from zero)
 since $(1 - \varkappa_n)$ is not greater than
$C\, |\chi^E(X_n)|$ and, by assumption 5), $C\, |\chi^E(X_n)|$ is
dominated by a decreasing geometrical progression.

Choose $\delta = \delta (\varepsilon, \varphi)$ such that:
$$
\omega_{\delta} (\varphi) :=
    \sup\limits_{\text{\rm dist}(x,y) < \delta} |\varphi(x) - \varphi(y)| <
    \varepsilon.
$$
By assumption 4), we have  $\sum_{n=1}^{\infty}\gamma_n<\infty$.
Choose $N = N (\varepsilon, \varphi)$ such that the following
holds:
$$
\sum\limits_N^{\infty} \gamma_k < \delta (\varepsilon, \varphi)
\qquad \mbox{and} \qquad \prod\limits_N^{\infty} \varkappa_k  > 1
- \varepsilon.
$$

Since the number of points in a pre-image for projections $\rho_n$
does not depend on a point in the image, a set ${\widetilde
X}_{m,\varepsilon} \subset X_m$
 where the total projection

$$
\rho_{m,N} = \rho_{m-1} \circ \dots \circ \rho_N\colon {\widetilde
X}_{m,\varepsilon}\to X_N
$$
is defined
 contains most of the orbit $X_m$:
$$
\frac { \# {\widetilde X}_{m,\varepsilon} } {\# X_m} \ge
    \prod \limits_{k=N}^{m-1} \varkappa_k \geq
    \prod\limits_N^{\infty} \varkappa_k >
    1 - \varepsilon.
$$
This implies 1) in Lemma~\ref{ergodicitylemma}.

Take arbitrary $m$ and $n$ with $m > n > N (\varepsilon,
\varphi)$. By construction, on the set ${\widetilde
X}_{m,\varepsilon}$ the total projection $\rho_{m,n} = \rho_{m-1}
\circ \dots \circ \rho_n$ is defined and
 for every point $x$ from the set
 $\widetilde{X}_{m,\varepsilon}\subset X_m$ we have
 $$
\text{\rm dist} (f^j(x), f^j(\rho_{m,\,n}(x)))<\delta(\varepsilon,
\varphi), \qquad \text{for all $j=0, 1, \ldots, \pi(X_n)-1$.}
$$
 Hence for $x \in {\widetilde X}_{m,\varepsilon}$ we have:
$$
\left| \int\varphi \, d\nu_{\pi(X_n)}(x) - \int\varphi \, d\mu_n
\right| <
    \omega_{\delta} (\varphi) <
    \varepsilon.
$$
Thus all conditions of Lemma~\ref{ergodicitylemma} are verified.
Therefore every limit point $\mu$ of  the sequence $\{\mu_n\}$ is
ergodic. By assumption 5), $\chi^E(\mu_n)\to 0$ as $n\to \infty$,
and, by Lemma~\ref{l.limit}, we have $\chi^E(\mu)=0$, that is,
$\mu$ is an ergodic invariant non-hyperbolic measure of $f$. To
prove Proposition~\ref{p.firstmainproposition}
 it remains to check that the support of $\mu$ is uncountable. We need
the following lemma. Denote by $U_\delta(x)$ the ball of radius
$\delta$ centered at $x$.

\blm \label{l.mpositive} Set $r_n=\sum\limits_{k=n}^{\infty}
  \gamma_k $. For every point $x\in X_n$ 
one has $\mu(\overline{U_{r_n}(x)})>0$. \elm

We postpone the proof of Lemma \ref{l.mpositive} and complete the
proof of the proposition assuming that Lemma \ref{l.mpositive}
holds. By assumption 4),
 we have
$$
r_n=\sum\limits_{k=n}^{\infty}
  \gamma_k <
\sum\limits_{k=n}^{\infty} \frac{\min_{j\le k} d_j}{3\cdot 2^k}\le
\frac{d_n}3.
$$
Thus, by the choice of $d_n$, any two closed $r_n$--balls with
centers at different points of an
  orbit~$X_n$ are disjoint. But, by Lemma~\ref{l.mpositive}, $\mu$-measure of each of these balls is positive, hence
  measure $\mu$
  can not be supported on less than $\pi(X_n)$ points. On the other hand, $n$ is arbitrary, and periods $\pi(X_n)$
  tend to infinity.
  Therefore measure $\mu$  can not be supported on a finite set.
  Since an infinite support of an invariant ergodic non-atomic measure is a
  closed set without isolated points, it can not be countable.
  This completes the proof of Proposition~\ref{p.firstmainproposition}.

\begin{proof}[Proof of Lemma~\ref{l.mpositive}]
Take any $n\in \NN$ and any point $x \in X_n$. Note
 that in its $\gamma_{n}$--neighborhood ($\gamma_n$ as in condition 3)) there are at least
$\frac{\# {\widetilde X}_{n+1, \varepsilon}}{\pi(X_n)}$ points of
the orbit $X_{n+1}$, where
$$
\frac{\# {\widetilde X}_{n+1, \varepsilon}}{\pi(X_{n})} = \frac{\# {\widetilde X}_{n+1, \varepsilon}}{\pi(
X_{n+1})} \,\, \frac{\pi ( {X}_{n+1})}{ \pi( X_{n})} \ge \varkappa_{n}\,
 \frac{\pi(X_{n+1})}{\pi(X_n)}=\bar \varkappa_{n}.
$$
Therefore, for all $n\in \NN$ and every point $x \in X_n$,
\begin{equation}
\label{e.medida}
 \mu_{n+1}(U_{\gamma_{n}}(x))\ge
\frac{\varkappa_{n}\,
 \frac{\pi(X_{n+1})}{\pi(X_n)}}{\pi(X_{n+1})}=
\varkappa_{n}\, \frac{1}{\pi(X_n)}= \varkappa_{n}\, \mu_n (\{x\}).
\end{equation}
Notice that in the neighborhood $U_{\gamma_n}(x)$ there are $p$
different points $x_1,\dots, x_p$ of the orbit $X_{n+1}$, where
$p\ge \bar\varkappa_n$. Notice also that by the definition of
$\ga_{n+1}$  the family of neighborhoods
$\{U_{\gamma_{n+1}}(x_i)\}_{i=1}^p$ is pairwise disjoint and their
union is contained in $U_{\gamma_{n+1}+\gamma_n} (x)$. Therefore,
since $p\ge \bar\varkappa_n$, and assuming that $x_{j_0}$
minimizes the measure $\{\mu_{n+2}
(U_{\gamma_{n+2}}(x_i))\}_{i=1}^p$, we have
$$
\mu_{n+2}(U_{\gamma_{n+1}+\gamma_n} (x)) \ge \sum_{i=1}^p
\mu_{n+2} (U_{\gamma_{n+1}}(x_i))\ge \left( \varkappa_{n}\,
 \frac{\pi(X_{n+1})}{\pi(X_n)}\right) \,\mu_{n+2}(U_{\gamma_{n+2}}(x_{j_0})).
$$
Since $x_{j_0}\in X_{n+1}$, equation~(\ref{e.medida}) now gives
$$
\begin{array}{ll}
\mu_{n+2}(U_{\gamma_{n+1}+\gamma_n} (x)) &\ge \left(
\varkappa_{n}\,
 \dfrac{\pi(X_{n+1})}{\pi(X_n)}\right)\, \varkappa_{n+1}\,
 \mu_{n+1}(x_{j_0})=\\\\ &=
 \left(
\varkappa_{n}\,
\dfrac{\pi(X_{n+1})}{\pi(X_n)}\right)\,\varkappa_{n+1}\,
 \dfrac{1}{\pi(X_{n+1})}=\displaystyle{\dfrac{\varkappa_n\,
 \varkappa_{n+1}}{\pi(X_{n})}}=\\\\&=
\varkappa_n\,
 \varkappa_{n+1}\, \mu_n(\{x\}).
 \end{array}
$$
Thus arguing inductively we have, for every $x\in X_n$,
$$
\mu_{n+\ell}(U_{\gamma_{n+\ell}+\cdots+
\gamma_{n+1}+\gamma_{n}}(x)) \ge (\varkappa_{n+\ell}\, \cdots\,
\varkappa_{n+1}\, \varkappa_{n}) \, \mu_n (\{x\}).
$$
Taking a limit and recalling equation (\ref{e.varkappa}), we have:
$$
\mu(\overline{U_{r_n} (x)}) \ge \left( \prod_{k=n}^{\infty}
\varkappa_{k}\right)\,
 \mu_n(\{x\}) >0, \quad \text{where $r_n=\sum\limits_{k=n}^{\infty}
\gamma_k$.} $$ Therefore, Lemma~\ref{l.mpositive} holds.
\end{proof} The proof of Proposition~\ref{p.firstmainproposition}
is now complete.
  \end{proof}



\section{Homoclinic classes of $C^1$-generic diffeomorphisms} \label{s.generic}

Here we state some properties of homoclinic classes of $C^1$-generic diffeomorphisms that
we will use  
 to obtain periodic points having Lyapunov exponents
close to zero. For this a key step is to get heterodimensional
cycles associated to saddles in these non-hyperbolic  homoclinic
class.


\subsection{Generic properties}
\label{ss.generic}

There is a residual subset ${\cal{G}}$ of $\diffM$ such that every
diffeomorphism $f\in \cG$ satisfies properties {\bf R1)}--{\bf
R4)}  below.
\begin{description}
\item{\bf R1)}
Every homoclinic class of $f\in \cG$ containing saddles of
indices $a$ and
$b$, $a<b$, also contains saddles of index $c$, for every
$c\in (a,b)\cap\NN$. See \cite[Theorem 1]{ABCDW}.
\end{description}

Consider a hyperbolic periodic point $P_f$ of a diffeomorphism
$f$. It is well known that there are open neighborhoods $U$ of
$P_f$ in the  manifold and $\mathcal{U}$ of $f$ in $\diffM$ such
that every $g\in \cU$ has a unique hyperbolic periodic point $P_g$
of the same period as $P_f$ in $U$. The point $P_g$ is called the
{\it continuation} of $P_f$.

\begin{description}
\item{\bf R2)}
Given any $f\in \cG$ and any pair of saddles $P_f$ and $Q_f$ of
$f$, there is a neighborhood $\cU_f$ of $f$ such that either
$H(P_g,g)=H(Q_g,g)$ for all $g\in \cU_f\cap \cG$, or $H(P_g,g)\cap
H(Q_g,g)=\emptyset$ for all $g\in \cU_f\cap \cG$. See \cite[Lemma
2.1]{ABCDW} (this lemma follows from \cite{CMP,BC} using a
 standard genericity argument). If the first case
holds, $H(P_g,g)=H(Q_g,g)$ for all $g\in \cU_f\cap \cG$, we say that the saddles $P_f$ and $Q_f$ are
\emph{persistently linked in $\cU_f$.}
\end{description}

Clearly, homoclinically related saddles are persistently linked.
 In fact, to be homoclinically related is
stronger than to be persistently linked. Homoclinically related
saddles have the same index while there are persistently linked
saddles having different indices (thus these saddles are not
homoclinically related).

\medskip

In order to state the last two generic conditions  we need some
general facts about homoclinic classes of hyperbolic points with
real multipliers.

\bdef Let $P$ be a periodic point of period $\pi(P)$ of a 
 diffeomorphism $f$. We say that
$P$ has {\emph{real multipliers}} if every eigenvalue $\lambda$ of
$Df^{\pi(P)}(P)$ is real and has multiplicity one, and two
different eigenvalues of $Df^{\pi(P)}(P)$ have different absolute values.
We order the eigenvalues of $Df^{\pi(P)}(P)$ in increasing
ordering according their absolute values $|\la_1(P)|<\cdots
<|\la_{\mathbf{m}}(P)|$ and say that $\la_k(P)$ is the
\emph{$k$-th multiplier of $P$}.
\endef

Consider a saddle $P$ with real multipliers with $s+1$ contracting
eigenvalues. Consider the bundle $E^{\sst}\subset T_PM$
corresponding to the first $s$ contracting eigenvalues of $P$ and
the {\emph{strong stable manifold}} $W^{\sst}(P)$ of $P$ (defined
as the only invariant manifold of dimension $s$ tangent to the
strong stable direction $E^{\sst}$).

\bdef\label{biacc}We say that a periodic point  $P$ of saddle type
with real multipliers  is {\emph{$\st$-biaccumulated}} (by its
homoclinic points) if both connected components of
$W^{\st}_{\loc}(P)\setminus W^{\sst}_{\loc}(P)$ contain transverse
homoclinic points of $P$. Define also {\emph{$\ut$-biaccumulation
by homoclinic points}} in a similar way.
\endef
\brm  Notice that $\st$- and $\ut$-biaccumulation are open
properties. \label{r.biopen}
 \erm

Note that if the homoclinic class of a saddle $P$ is non-trivial
then there is some  transverse homoclinic point $Z$ associated to
$P$. By the Smale homoclinic theorem, there is a small
neighborhood $U$ of the orbits of $P$ and $Z$ such that the
maximal invariant set $\La_f(U)$ of $f$ in   $U$ is a ``horseshoe"
(non-trivial hyperbolic set). Moreover, if $P$ has real
multipliers,  we can assume (after a perturbation if needed) that
the there is a $Df$-invariant splitting of $T_{\La_f(U)} M$ of the
form $E^{\sst}\oplus E^{\cst} \oplus E^{\ut}$, where $E^{\sst}$ is
a strong stable bundle, $E^{\cst}$ is a one-dimensional stable
bundle, and $E^{\ut}$ is a unstable bundle. Then for any periodic
point $A\in \La_f(U)$ there is defined its strong stable manifold
$W^{\sst}_{\loc}(A)$, thus the notion of $\st$-biaccumulation is
well-defined for every periodic point in $\La_f(U)$.

The density of periodic points homoclinically related to $P$ in
the set $\La_f(U)$ and the structure of local product immediately
imply the following result.

\blm \label{l.biaccumulated} Let $P$ be a hyperbolic periodic
point of $f$ having real multipliers whose homoclinic class is
nontrivial. Then there is a $C^1$-open set $\,\cU$, $f$ is in the
closure of $\,\cU$, such that
every $g\in \cU$ has hyperbolic saddles from $\La_g(U)$ which are
homoclinically related to $P_g$ and are $\st$-biaccumulated by
their transverse homoclinic points. Similarly, for
$\ut$-biaccumulation. \elm

 We can now formulate the last two genericity conditions that we need.
Given a saddle $P$, we denote by $\mbox{Per}_{\RR}(H(P,f))$ the
saddles homoclinically related to $P$ having real multipliers.
 Clearly,
 $\mbox{Per}_{\RR}(H(P,f))\subset H(P,f)$.

\begin{description}
\item{\bf R3)}
For every diffeomorphism $f\in \cG$ and every saddle $P$ of $f$
whose homoclinic class is non-trivial the set
$\mbox{Per}_{\RR}(H(P,f))$ is dense in the whole homoclinic class
$H(P,f)$. See \cite[Proposition 2.3]{ABCDW}, which is just a
dynamical reformulation of \cite[Lemmas 1.9 and 4.6]{BDP}.
Moreover, by Lemma~\ref{l.biaccumulated} we can also assume that
there is a dense set of  points in $\mbox{Per}_{\RR}(H(P,f))$
which are $\st$- (or $\ut$-) biaccumulated.
\item{\bf R4)}
Consider an open set $\cU$ of $\diffM$ such that there are
hyperbolic saddles $P_f$ and $Q_f$ which are persistently linked
in $\cU$. Suppose that the dimensions of their stable bundles are
$s+1$ and $s$, respectively. Then for every $f$ from the residual
subset $\cG\cap \cU$ of $\cU$ and for every $\varepsilon>0$ the
sets
$$
\begin{array}{ll}
\mbox{Per}^{(1-\varepsilon, 1)}_{\RR}(H(P_f,f))&=\{R_f\in \mbox{Per}_{\RR}(H(P_f,f)) \ : \ |\la_{s+1}(R_f)|\in (1-\varepsilon,1)\},\\
\mbox{Per}^{(1, 1+\varepsilon)}_{\RR}(H(Q_f,f))&=\{R_f\in \mbox{Per}_{\RR}(H(Q_f,f)) \ : \ |\la_{s+1}(R_f)|\in (1,1+\varepsilon)\}
\end{array}
$$  are dense in
$\mbox{Per}_{\RR}(H(P_f,f))$ and $\mbox{Per}_{\RR}(H(Q_f,f))$,
respectively. This condition is
 an immediate consequence
of the results in \cite{BDF}.
 By Lemma \ref{l.biaccumulated}, these saddles $\{R_f\}$ can be taken also with the
biaccumulation property.
\end{description}



\subsection{Creation of cycles}
\label{ss.creationofcycles}

In this section, we state results that allow us to generate
heterodimesional cycles for persistently linked saddles.

\bpro Let $\, \cU$ be an open set of $\,\diffM$ such that there
are saddles $P_f$ and $Q_f$ (depending continuously on $f\in \cU$)
with consecutive indices  which are persistently linked in \ $\cU$.
Then there is a dense subset $\cH$ of $\,\cU$ such that every
diffeomorphism $f\in \cH$ has a coindex one heterodimensional
cycle associated to saddles $A_f$ and $B_f$ such that
\begin{itemize}
\item
the saddles $A_f$ and $B_f$ have real multipliers,
\item
the saddle $A_f$ is homoclinically related to $P_f$ and the saddle
$B_f$ is homoclinically related to $Q_f$.
\end{itemize}
\label{p.bdhayashi} \epro

\begin{proof}
First, note that since the saddles $P_f$ and $Q_f$ are
persistently linked, their homoclinic classes are both
non-trivial.  Note also that it is enough to
prove this result in a small neighborhood $\cV$ of $f\in \cU\cap
\cG$ ($\cG$ is the residual set of $\diff (M)$ in
Section~\ref{ss.generic}). By condition {\bf R3)}, every
diffeomorphism $f\in \cG\cap\cU$ has a pair of saddles $A_f$ and
$B_f$ with real multipliers and which are homoclinically related
to $P_f$ and $Q_f$, respectively. In particular, the saddles $A_f$
and $B_f$ verify $H(A_f,f)=H(P_f,f)$ and $H(B_f,f)=H(Q_f,f)$.
Thus, by the definition of persistently linked saddles, after
shrinking $\cV$ we can assume that
$$
H(A_f,f)=H(P_f,f)=H(Q_f,f)=H(B_f,f), \quad \mbox{for every $f\in
\cV\cap \cG$.}
$$

We now  get  heterodimensional cycles associated to $A_f$ and
$B_f$. We use a standard argument which follows by applying twice
Hayashi's Connecting Lemma to the saddles $A_f$ and $B_f$:

\blm [Hayashi's Connecting Lemma, \cite{H}] Consider a
diffeomorphism $f$ with a pair of saddles $M_f$ and $N_f$ such
that there are sequences of points $T_i$ and of natural numbers
$n_i$ such that $T_i$ accumulates to $W^{\ut}_{\loc}(M_f,f)$ and
$f^{n_i}(T_i)$ accumulates to $W^{\st}_{\loc}(N_f,f)$.
 Then there is $g$ arbitrarily $C^1$-close to $f$ such that
$W^{\ut}(M_g,g)\cap W^{\st}(N_g,g)\ne \emptyset$.
\label{l.hayashi}
 \elm

Note that this lemma can be applied to any pair of saddles $M_f$
and $N_f$ in the same transitive set of $f$ (for instance, a
homoclinic class).

We observe that our arguments are now local, thus by shrinking
 $\cV$, we can assume that the saddles
$A_f$ and $B_f$ are defined in the whole $\cV$. Consider the
subsets of $\cV$ defined by
$$
\begin{array}{ll}
\cI&=\{g\in \cV \colon W^{\st}(A_g,g)\cap
W^{\ut}(B_g,g)\ne\emptyset\} \quad
\mbox{and}\\
\cJ&=\{g\in \cV \colon W^{\ut}(A_g,g)\cap
W^{\st}(B_g,g)\ne\emptyset\}.
\end{array}
$$
Since the set $H(A_g,g)=H(B_g,g)$, $g\in \cG$, is transitive, we
can apply Lemma~\ref{l.hayashi} to the saddles $A_g$ and $B_g$,
obtaining that the sets $\cI$ and $\cJ$  are both dense in $\cV$.

Suppose for a moment that the index of $A_g$ is less than the
index of $B_g$. Then the set $\cI$ has non-empty interior and
therefore $ \cI \cap \cJ$ is dense in $\cV$. Indeed, notice that
the sum of the dimensions of the stable manifold of $A_g$ and the
unstable manifold of $B_g$ is greater than the dimension of the
ambient manifold. Thus, after a perturbation, an intersection
between $W^{\st}(A_g,g)$ and $W^{\ut}(B_g,g)$ can be made a
transverse one, thus persistent after perturbations. Hence the set
$\cI$ contains an open an dense subset $\cY$  of $\cV$.

Finally, consider the set $\cH=\cY\cap \cJ\subset \cI\cap \cJ$. By
the previous construction, $\cH$ is dense in $\cV$. Finally, by
definition of $\cI$ and $\cJ$, the set $\cH$ consists of
diffeomorphisms $g$ with a heterodimensional cycle associated to
the saddles $A_g$ and $B_g$. The proof of the proposition is now
complete.
\end{proof}


\section{Generation of cycles and saddles with central exponents
close to zero} \label{s.generationzero}

In this section, we state two  technical propositions about
generation of (heterodimensional) cycles and of saddles with
Lyapunov exponents close to zero inside non-hyperbolic homoclinic
classes of generic diffeomorphisms. In Section~\ref{s.genzero}, we
will combine these results and the ones in
Section~\ref{s.ergodicmeasures} to get non-hyperbolic ergodic
measures with uncountable support.

Below we will restrict our attention to the dynamics in an open
set $V$ of $M$. Recall that $\La_f(V)$ is the maximal invariant
set of $f$ in $V$, $\La_f(V)=\cap_{i\in \mathbb{Z}} f^i(V)$.

\bdef[$V$-relative dynamics]$\,$
\begin{itemize}
 \item
A pair of saddles $A$ and $B$ of different indices have a
{\emph{$V$-related (heterodimensional) cycle}} if the set $V$
contains the orbits of $A$ and $B$ and there are heteroclinic
points $X\in W^{\ut}(A)\cap W^{\st}(B)$ and $Y\in W^{\st}(A)\cap
W^{\ut}(B)$ whose orbits are contained in $V$.
\item
Two saddles $A$ and $A^\prime$ are {\emph{$V$-homoclinically
related}} if $V$ is a neighborhood of the orbits of $A$ and
$A^\prime$ and there are heteroclinic  points  $Z\in
W^{\st}(A)\cap W^{\ut}(A^\prime)$ and $Z^\prime\in W^{\ut}(A)\cap
W^{\st}(A^\prime)$ whose orbits are contained in $V$.
\item
Consider a saddle  $A$ with real multipliers
and a neighborhood $V$ of it. The saddle $A$ is
 {\emph{$V$-$\st$-biaccumulated}}
if both connected components of $W^{\st}_{\loc}(A)\setminus
W^{{\sst}}_{\loc}(A)$ contain homoclinic points of $A$ whose orbits
are contained in $V$. We define {\emph{$V$-$\ut$-biaccumulation}}
similarly.
\end{itemize}
\endef

In order to state the two main technical results of this section,
we need the following lemma (which easily follows
from the $\la$-lemma and the Smale's homoclinic theorem; we present the proof at the end of Section~\ref{ss.simple})
which allows us to identify a fixed region of the manifold
where the dynamics in a neighborhood of a cycle is partially hyperbolic:



\blm\label{l.partialhyperbolicity} Let $f$ be a diffeomorphism with
a heterodimensional 
 cycle associated to saddles $A_f$ and $B_f$
such that
\begin{itemize}
\item the saddles $A_f$ and $B_f$ have real multipliers
 and \ $\mbox{\rm index\,}(A_f)+1=\mbox{\rm
index\,}(B_f)$;
\item
$A_f$ is $\st$-biaccumulated and $B_f$ is $\ut$-biaccumulated.
\end{itemize}
Then arbitrarily $C^1$-close to $f$ there are diffeomorphisms $g$
such that for some open set $V\subseteq M$
\begin{itemize}
\item
the set $\Lambda_g({V})$ is partially hyperbolic (with a splitting
$E^{\sst}\oplus E^{\cent} \oplus E^{\uut}$, $E^{\cent}$ is
one-dimensional);
\item
$A_g$ is ${V}$-$\st$-biaccumulated and $B_g$ is
${V}$-$\ut$-biaccumulated;
\item
the diffeomorphism $g$ has a ${V}$-related cycle associated to $A_g$
and $B_g$.
\end{itemize}
 \elm

Lemma~\ref{l.partialhyperbolicity} will be used in
Section~\ref{s.genzero} just once, at the beginning of the
construction, to get an open set where the relative dynamics is
partially hyperbolic. The  two propositions below will be used on
each step of the inductive construction in Section~\ref{s.genzero}.

\bpro\label{p.simpleanddense} Let $f$ have a $V$-related
 cycle associated to saddles $A_f$ and $B_f$ such that
\begin{itemize}
\item
$A_f$ and $B_f$ have real multipliers and $\mbox{\rm index\,}(A_f)+1=\mbox{\rm
index\,}(B_f)$;
\item
$A_f$ is $V$-$\st$-biaccumulated and $B_f$ is
$V$-$\ut$-biaccumulated;
\item
the set $\Lambda_g({V})$ is partially hyperbolic (with a splitting
$E^{\sst}\oplus E^{\cent} \oplus E^{\uut}$, $E^{\cent}$ is
one-dimensional).
\end{itemize}
 Then arbitrarily $C^1$-close to $f$ there
is an open set $\mathcal{E}\subset \text{\rm Diff}^{\,1}(M)$
exhibiting a dense subset $\mathcal{D}\subset \mathcal{E}$ such
that every diffeomorphism $g\in \mathcal{D}$ has a $V$-related cycle
associated with $A_g$ and $B_g$. \epro

\brm Weaker versions (without  biaccumulation hypotheses) of
Proposition~\ref{p.simpleanddense} can be obtained following
{\rm\cite{BD4}}. Here we need two extra conclusions which
do not follow  straightforwardly from {\rm\cite{BD4}}
and whose proofs demand some adaptations.
We first need that
 the persistent cycles were
associated to the continuation of the initial saddles.
 Second, these cycles must be
$V$-related. To get these conclusions we use
the biaccumulation hypotheses. 
\erm

Given an open set $V$,
we say that two invariant manifolds of a diffeomorphism are {\emph{$V$-related}}
if they have  an intersection point whose orbit is contained in $V$.
The following extension of \cite[Proposition 4.1]{BDF} is
the main technical step in the proof of Proposition~\ref{p.simpleanddense}. It will also be used in Section~\ref{s.genzero}.

\bpro \label{p.bdf} Let $f$ be a
diffeomorphism with a $V$-related  cycle
associated to saddles $A_f$ and $B_f$ such that
\begin{description}
\item{{\bf (i)}}
the saddles $A_f$ and $B_f$ have real multipliers;
\item{{\bf (ii)}}
$\text{\rm index}(A_f)+1=\text{\rm index}(B_f)$;
\item{{\bf (iii)}}
$A_f$ is $V$-$\st$-biaccumulated, and  $B_f$ is
$V$-$\ut$-biaccu\-mu\-la\-ted.
\end{description}

Then there are  sequences of natural numbers $\ell_k,m_k$, that tend
to infinity as $k\to \infty$, and a sequence of diffeomorphisms
$f_k$,
 $f_k\to f$ as $k\to \infty$, such that $f_k$ coincides with $f$ along the orbits of $A_{f}$ and $B_{f}$,
 and has a hyperbolic saddle $R_{k}\in \Lambda_{f_{k}}(V)$ having real multipliers  with the following properties:
\begin{description}
\item{{\bf (1)}}
fix  neighborhoods $U_B$ of the orbit of $B_f$ and   $U_A$ of the
orbit of $A_f$; the orbit of the saddle $R_{k}$ spends a fixed
number $t_{(a,b)}$ (independent of $k$) of iterates
  to go from  $U_B$ to $U_A$, then it remains $\ell_k\,\pi(A_f)$ iterates in  $U_A$, then it takes
a fixed number of iterates $t_{(b,a)}$ (independent of $k$) to go
from $U_A$ to $U_B$, and finally it remains $m_k\,\pi(B_f)$
iterates in $U_B$. In particular, there is a  constant $t\in \NN$
 independent of $k$ such that the period of $R_k$ is
$
\pi(R_{k})=m_k\, \pi(B_{f})+ \ell_k \, \pi(A_{f})+ t;
$
\item{{\bf (2)}}
there is a constant $\const>0$ independent of $k$ such that the
central multiplier of $R_{k}$\footnote{Let
$\ell=\mathbf{m}-\text{\rm index}(B_f)$, the central multiplier of
$R_k$ is its $\ell$-th multiplier, where $\mathbf{m}$ is the
dimension of the ambient manifold. See also Definition \ref{d.centralmult}.}
  satisfies $ {\const}^{-1}<|\la^{\ct}
(R_{k})|<\const.$ 
\end{description}
Suppose also that the central multiplier $\la^{\ct}(A_{f})$ of
$A_f$ is close to one. Then
 \begin{description}
\item{{\bf (3)}}
 $R_{k}$ has the same index as $B_{f}$ and is 
$V$-homoclinically related to $B_{f}$;
\item{{\bf (4)}}
$W^{\st}(R_k)$ and $W^{\uut}(R_k)$ are $V$-related, and
$W^{\uut}(R_k)$ and $W^{\st}(B_f)$ are $V$-related (these
intersections are
 quasi-transverse);
\item{{\bf (5)}}
there is a $V$-related cycle associated to $R_{k}$ and $A_{f}$.
\end{description}
\epro

\begin{figure}[htb]

\begin{center}
\psfrag{R}{$R_k$} \psfrag{W}{$W^{\ut}(R_k)$} \psfrag{A}{$A_f$}
\psfrag{B}{$B_f$}
\includegraphics[height=1.6in]{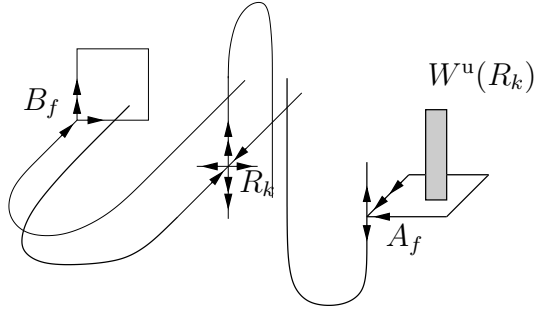}
 \caption{Dynamical configuration in Proposition~\ref{p.bdf}}
\label{f.configuration}
\end{center}

\end{figure}

The dynamical configuration of Proposition~\ref{p.bdf} is depicted in
Figure~\ref{f.configuration}.

We will prove Propositions~\ref{p.simpleanddense} and \ref{p.bdf} in Section~\ref{ss.proofs}.
In order to do that, we first consider special
cycles  (the so-called {\emph{simple cycles}})
and study their  unfolding by special parametrized families of
diffeomorphisms.  In Section~\ref{s.genzero},
using Propositions~\ref{p.bdf} and \ref{p.simpleanddense}, we will conclude the proof of Theorem~\ref{t.main}.

\subsection{Simple heterodimensional cycles}
\label{ss.simple} We adapt here some constructions from
\cite{ABCDW,BD4,BDF}. The details of these constructions can be
found in \cite[Section 3.1]{ABCDW} and
\cite[Section 3]{BD4}. 

\bdef[Co-index one cycles and central multipliers] Let  $f$ be a
diffeomorphism with a heterodimensional cycle associated to
saddles $A$ and $B$.
\begin{itemize}
\item
The cycle has \emph{co-index one} if
\,$\text{{\rm index\,}}(A)=\text{{\rm index\,}}(B)\pm 1$.
\item
Assume that the saddles $A$ and $B$ have real multipliers and let
$(s+1)$ and $s$ be the dimensions of the stable bundles of $A$ and
$B$,  respectively.
The {\em central multipliers of the cycle} are the $(s+1)$-th
multipliers of $A$ and $B$,  denoted  by $\la^{\ct}(A)$ and
$\la^{\ct}(B)$, respectively.
\end{itemize}
\label{d.centralmult}
\endef

Consider a diffeomorphism $f$ with a co-index one cycle associated to saddles
$A$ and $B$ with real multipliers as above.
Fix heteroclinic points
of the cycle
 $$
X\in W^{\st}(A,f)\cap W^{\ut}(B,f)\quad  \text{\rm and}\quad Y\in
W^{\ut}(A,f)\cap W^{\st}(B,f)
$$
and a small neighborhood $\neig$ of the orbits of the saddles $A$ and $B$ in the cycle and
the heteroclinic points $X$ and $Y$. In this way, we get a $\neig$-related cycle.

Notice that the saddles $A$ and $B$ have  indices   $(\mathbf{m}-s-1)$ and $(\mathbf{m}-s)$, where $\mathbf{m}$
is the dimension of the ambient manifold $M$. Since the saddles have real multipliers, there is a (unique)
$Df$-invariant dominated splitting defined on the union of the orbits  $\cO_{A}$ of $A$ and $\cO_{B}$ of $B$,
$$
T_x M=E^{\sst}_x\oplus E^{\ct}_x\oplus E^{\uut}_x, \qquad x\in
\cO_{A} \cup \cO_{B},
$$
where $\dim E^{\ct}_x=1$, $\dim E^{\sst}_x=s$, and $\dim
E^{\uut}_x=(\mathbf{m}-s-1)=u$. After a perturbation (while
keeping the cycle), we can assume that there are neighborhoods
$U_A$ and $U_B$ of $\cO_A$ and $\cO_B$, contained in the set
$\neig$, and coordinates in these neighborhoods where $f^{\pi(A)}$
and $f^{\pi(B)}$ are linear maps, and the splitting
$E^{\sst}\oplus E^{\ct}\oplus E^{\uut}$ is of the form
$$
E^{\sst}=\RR^s\times\{(0,0^u)\}, \quad E^{\ct}=\{0^s\}\times
\RR\times\{0^u\}, \quad E^{\uut}=\{(0^s,0)\}\times \RR^u.
$$

Observe that the sum of the dimensions of $W^{\st}(A,f)$ and
$W^{\ut}(B,f)$ is $(\mathbf{m}+1)$. Thus, after another
perturbation, we can assume that the intersection at the
heteroclinic point $X\in W^{\st}(A,f)\cap W^{\ut}(B,f)$ is
transverse. Similarly, we can also assume that the intersection
between $W^{\ut}(A,f)$ and $W^{\st}(B,f)$ at $Y$ is
quasi-transverse, i.e., $T_Y W^{\ut}(A,f)+T_Y W^{\st}(B,f)= T_Y
W^{\ut}(A,f)\oplus T_Y W^{\st}(B,f)$ and this sum has dimension
$(\mathbf{m}-1)$.

Take small neighborhoods $U_X$ and $U_Y$ of the heteroclinic
points $X$ and $Y$ and natural numbers $n$ and $m$ such that
 $$
 f^{n} (U_X)\subset U_A, \quad f^{-n} (U_X)\subset U_B, \quad
f^{-m} (U_Y)\subset U_A,
 \quad
 \mbox{and}
 \quad
\quad f^{m} (U_Y)\subset U_B
$$
and
$$
\bigcup_{i=-n}^n f^i(U_X)\subset \neig \quad \mbox{and}
 \quad
\bigcup_{i=-m}^m f^i(U_Y)\subset \neig.
$$

Consider the \emph{transition times} $t_{(b,a)}=2\,n$ and
$t_{(a,b)}=2\,m$ and define \emph{transition maps} $\mfT_{(a,b)}$
from $U_A$ to $U_B$ and $\mfT_{(b,a)}$ from $U_B$ to $U_A$ defined
on small neighborhoods $U_Y^-\subset U_A$ of $f^{-m}(Y)=Y^-$
 and  $U_X^-\subset U_B$ of $f^{-n}(X)=X^-$ as follows
$$
\mfT_{(a,b)}=f^{t_{(a,b)}}  \colon U_{Y}^-\to U_B
\quad \mbox{and} \quad
 \mfT_{(b,a)}=f^{t_{(b,a)}} \colon U_{X}^-\to U_A.
$$
After a perturbation, we can assume that the $\neig$-related cycle
is \emph{simple}\footnote{We use the notation in \cite[Definition
3.5]{BD4} corresponding to the affine cycles in \cite{ABCDW}. The
difference between these two definitions is that in \cite{BD4} the
central components $T^{\ct}_{(a,b)}$ and $T^{\ct}_{(b,a)}$ of the
transitions are isometries while in \cite{ABCDW} are just affine
maps.}.  This means that in the local coordinates in $U_A$ and
$U_B$ above one can write
$$
\mfT_{(b,a)}=
(T^{\st}_{(b,a)},T^{\ct}_{(b,a)},T^{\ut}_{(b,a)})\quad \mbox{and}
\quad \mfT_{(a,b)}=
(T^{\st}_{(a,b)},T^{\ct}_{(a,b)},T^{\ut}_{(a,b)})
$$
where
\begin{description}
\item{\bf S1)}
$T^{\st}_{(i,j)}\colon \RR^s\to \RR^s$ and
$(T^{\ut}_{(i,j)})^{-1}\colon \RR^u\to \RR^u$ are affine
contractions;
\item{\bf S2)}
$T^{\ct}_{(a,b)}\colon \RR\to \RR$ and $T^{\ct}_{(b,a)}\colon
\RR\to \RR$ are affine isometries. Moreover, $T^{\ct}_{(a,b)}$
 is linear (note that the central coordinates of the
heteroclinic points $f^{-m}(Y)$ and $f^m(Y)$ are both zero).
\end{description}

\begin{figure}[htb]

\begin{center}
\psfrag{Y}{$Y$}
\psfrag{Xp}{$X^+$}
\psfrag{Yp}{$Y^-$}
\psfrag{Yq}{$Y^+$}
\psfrag{Xq}{$X^-$}

\psfrag{Uyp}{$U_{Y}^-$}
\psfrag{Uxq}{$U_{X}^-$}
\psfrag{Q}{$B$}
\psfrag{P}{$A$}
\psfrag{T1}{$\mfT_{(a,b)}$}
\psfrag{T2}{$\mfT_{(b,a)}$}
\includegraphics[height=1.6in]{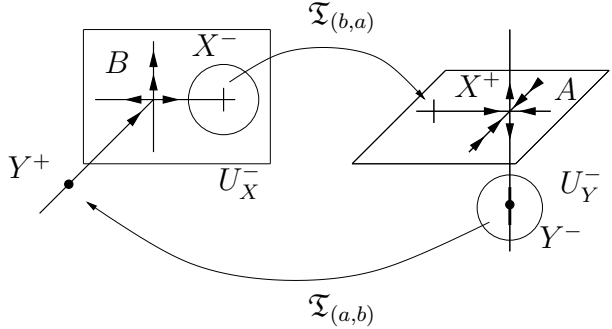}
 \caption{Simple cycle}
\label{f.simples}
\end{center}

\end{figure}

As the splitting $E^{\sst}\oplus E^{\ct}\oplus E^{\uut}$ is
 dominated,
properties {\bf S1)} and {\bf S2)} can be obtained making perturbations along the finite orbits $\{f^{-n}(X),\ldots,  f^{n}(X)\}$
  and $\{f^{-m}(Y),\ldots,
 f^{m}(Y)\}$, after increasing the transition times $n$ and $m$ and
shrinking the neighborhood $\neig$, if necessary.

The points $X\in W^{\st}(A,f)\cap W^{\ut}(B,f)$ and $Y\in
W^{\ut}(A,f)\cap W^{\st}(B,f)$ are the \emph{transverse} and
\emph{quasi-transverse heteroclinic intersections of the simple
cycle} and the set
$\neig$ is a
\emph{neighborhood of the simple cycle}. Note that the cycle associated
to the saddles $A$ and $B$ with heteroclinic orbits $X$ and $Y$ is
$\neig$-related. We refer to this sort of cycles as {\emph{$\neig$-related simple cycles.}}

The previous construction gives the (affine) dynamics of $f$ in
the neighborhood $\neig$ of the cycle. Note that he dynamics on
the maximal invariant set of $f$  in $\neig$ is partially
hyperbolic (with a splitting $E^{\sst}\oplus E^{\ct}\oplus
E^{\uut}$). We say that  $\neig$ is a partially hyperbolic
neighborhood. This property persists after perturbations.

Next lemma summarizes the construction above:

\blm{\bf (\cite[Lemma 3.4]{ABCDW}, \cite[Proposition 3.6]{BD4})}
Suppose  a diffeomorphism $f$ has a co-index one cycle associated to saddles $A_f$ and $B_f$ with real multipliers. Then there are an open set $\neig\subseteq M$ and
a diffeomorphism $g$ arbitrarily
$C^1$-close to $f$ having a  $\neig$-related simple cycle associated to $A_g$ and
$B_g$. Moreover,
the dynamics of $g$ in $\neig$ is
partially hyperbolic.
\label{l.aproxsimple}
\elm




\subsubsection{Sketch of the proof of Lemma \ref{l.partialhyperbolicity}}
This lemma is consequence of the previous constructions. By
Lemma~\ref{l.aproxsimple}, there is $g$ arbitrarily close to $f$ and
an open set $V_0\subseteq M$ such that $g$ has a
$V_0$-related simple cycle. In particular, the dynamics of $g$
in $V_0$ is partially hyperbolic.

Consider the saddle $A_g$ with real multipliers. Since the initial
saddle $A_f$ is
$\st$-biaccu\-mu\-la\-ted, 
there are transverse homoclinic points of $A_g$ in both components
of $(W^{\st}_{\loc}(A_g)\setminus W^{\sst}_{\loc}(A_g))$, recall
Remark~\ref{r.biopen}.  We
take now a small neighborhood $V_A$ of the orbits of these
transverse homoclinic points and of the orbit of $A_g$. By the
homoclinic theorem, if this neighborhood is small enough, the
maximal invariant set of $g$ in $V_A$ is hyperbolic,
and (after a
perturbation, if necessary) each bundle of this splitting is the
sum of one-dimensional invariant bundles.
 We define
the set $V_B$ similarly (using that $B_g$ is $\ut$-biaccumulated).
Now it suffices to consider ${V}=V_0 \cup V_A \cup V_B$.
\hfill$\Box$

\subsection{Unfolding of simple
cycles}\label{ss.simpleunfolding}

In this section, we  consider the unfolding of simple cycles via
special families of diffeomorphisms preserving the partially
hyperbolic structure of the  cycle. Our goal is to
generate  periodic
points with bounded central multipliers  (as the points $R_k$
in Proposition~\ref{p.bdf}). Later, in Section~\ref{ss.intersections},
we obtain these saddles satisfying some additional (homoclinic/heteroclinic)
intersection properties.

Consider a diffeomorphism $f$ having
a simple cycle associated to saddles $A$ and $B$ (since in this section these points
remain fixed we will omit the subscript denoting the dependence on the diffeomorphism),
$\mbox{index}(A)=\mbox{index}(B)-1$,  heteroclinic points $X$ (transverse) and $Y$ (quasi-transverse),
and associated neighborhood $\neig$, as in Section~\ref{ss.simple}.
 Following the notation in Section~\ref{ss.simple}, consider
the one-parameter family of transitions $(\mfT_{(a,b),\nu})_{\nu}$
from $A$ to $B$ that in the local coordinates has the form
$$
\mfT_{(a,b),\nu}=\mfT_{(a,b)}+(0^s,\nu,0^u).
$$
For every small $\nu$, there is  a perturbation $f_{\nu}$ of
$f$ in a vicinity of the heteroclinic point
$f^m(Y)=Y^+$ such that (in local coordinates) one has
$$
f_\nu^{t_{(a,b)}}(x^s,x,x^u)=
f^{t_{(a,b)}}(x^s,x,x^u)+(0^s,\nu,0^u).
$$
Thus the orbits of $A$ and $B$ are not modified and $A$ and $B$
are periodic points of $f_{\nu}$. Note also that $f_\nu^{t_{(b,a)}}$
coincides with $\mfT_{(b,a)}$ in the neighborhood $U_X^-$ of
$f^{-n}(X)=X^-$ (recall that $X$ is the transverse heteroclinic
point of the cycle).

 Note that by construction, for
all large $\ell$ and $m$, there is a small neighborhood
$U_{\ell,m}\subset U_B$ close to $X^-$ with
$$
(f^{\pi_{\ell,m}}_{\nu})_{|_{U_{\ell,m}}} = f^{m\, \pi(B)}\circ
\mfT_{(a,b),\nu} \circ f^{\ell\, \pi(A)}\circ \mfT_{(b,a)}\colon U_{\ell,m}\to U_B,
$$
where
$$
\pi_{\ell,m}=    m\, \pi(B) +t_{(a,b)}+  \ell\, \pi(A)+ t_{(b,a)}.
$$

\blm \label{l.Rlm}
For all large  $\ell$ and $m$, there is $\nu_{\ell,m}$,
$\nu_{\ell,m}\to 0$ as $\ell,m\to \infty$, such that
$f_{\nu_{\ell,m}}$ has a saddle $R_{\ell,m}\in
\La_{f_{\nu_{\ell,m}}}(\neig)$ such that in the coordinates  in $U_B$
one has $R_{\ell,m}=(r^s_{\ell,m},c,r^u_{\ell,m})\in U_B$,
$|r^{s,u}_{\ell,m}|\to 0$ as $\ell,m\to \infty$. The period
of $R_{\ell,m}$ is
$$
\pi(R_{\ell,m})=m\, \pi(B)+ \ell \, \pi(A)+ t_{(a,b)}+t_{(b,a)}
$$
and its central multiplier satisfies
$$
|\la^{\ct}(R_{\ell,m})|= |(\la^{\ct}(B))^m\, (\la^{\ct}(A))^\ell|.
$$
Moreover, if $R_{\ell,m}$ is hyperbolic then $\mbox{\rm index
}(R_{\ell,m})\in\{ \mbox{\rm index }(A), \mbox{\rm index }(B)\}$,
according to the absolute value of $\la^{\ct}(R_{\ell,m})$.
 \elm

\brm\label{r.item1} Due to the construction above the saddle $R_{\ell,m}$ satisfies
item 1)
 in Proposition~\ref{p.bdf}. To prove 2)--5) in
Proposition~\ref{p.bdf}, we will need a more accurate control of
the central dynamics of the cycle and use the biaccumulation
properties. For that we need to take a sequence $(\ell_k,m_k)$ of
the integers above and let $R_k=R_{\ell_k,m_k}$. \erm

\brm \label{r.quotient} By construction, in the neighborhood $\neig$
of the simple cycle the diffeomorphisms $f_\nu$ keep invariant the
codimension one foliation generated by the sum of the strong
stable and strong unstable directions (hyperplanes parallel to
$\RR^s \times  \{0\} \times \RR^u$). Moreover, every $f_\nu$ acts
hyperbolically on these hyperplanes. We  consider the quotient
dynamics by these hyperplanes. Periodic points of this quotient
(one dimensional)
dynamics correspond to periodic points of the diffeomorphism
$f_\nu$. If the periodic point of the quotient dynamics is
expanding the corresponding periodic point of $f_\nu$ has the same index as $B$.
\erm

\medskip
{\noindent{{\emph{Proof of Lemma~\ref{l.Rlm}:}}} Suppose for
simplicity that in our local coordinates
$$
X^-=f^{-n}(X)=(0^s,1,0^u)\in U_B, \quad
 X^+=f^{n}(X)=(0^s,-1,0^u)\in U_A.
$$
Note that $T^{\ct}_{(a,b)}(x)=\tau\,x$, where $\tau=\pm 1$. Fix
large $\ell$ and $m$ and $B$, and let
$$
\nu_{\ell,m}=(\la^{\ct}(B))^{-m}+  \tau\, (\la^{\ct}(A))^\ell,
\qquad \nu_{\ell,m}\to 0 \quad \mbox{as}\quad  \ell,m\to \infty.
$$
Therefore, by definition of $\nu_{\ell,m}$,
\begin{equation}
\label{e.beta}
 (\la^{\ct}(B))^m\,(-\tau \, (\la^{\ct}(A))^\ell +\nu_{\ell,m})=1.
\end{equation}
This choice and equalities
$$
T^{\ct}_{(b,a)}(1)=-1, \qquad
T^{\ct}_{(a,b),\nu_{\ell,m}}(x)=T^{\ct}_{(a,b)}(x)+\nu_{\ell,m}=\tau\,x
+\nu_{\ell,m}
$$
imply that $1$ is a fixed point of the quotient dynamics:
$$
\begin{array}{ll}
&(\la^{\ct}(B))^m \circ T^{\ct}_{(a,b),\nu_{\ell,m}} \circ
(\la^{\ct}(A))^\ell \circ T^{\ct}_{(b,a)}(1)=\\
&\qquad \qquad \qquad = (\la^{\ct}(B))^m \circ
T^{\ct}_{(a,b),\nu_{\ell,m}}
\circ (\la^{\ct}(A))^\ell (-1)=\\
&\qquad \qquad \qquad = (\la^{\ct}(B))^m \circ
T^{\ct}_{(a,b),\nu_{\ell,m}} (-(\la^{\ct}(A))^\ell)=
\\
&\qquad \qquad \qquad= (\la^{\ct}(B))^m\, (-\tau\,
(\la^{\ct}(A))^\ell +\nu_{\ell,m})=1.
\end{array}
$$
Since $f_{\ell,m}=f_{\nu_{\ell,m}}$ preserves the $E^{\sst}$,
$E^{\uut}$, and $E^{\ct}$ directions, the hyperbolicity of the
directions $E^{\sst}$ and $E^{\uut}$ implies that the map
$$
f^{m\, \pi(B)}\circ \mfT_{(a,b),\nu_{\ell,m}} \circ f^{\ell\,
\pi(A)}\circ \mfT_{(b,a)}
$$
has a fixed point $R_{\ell,m}=(r^s_{\ell,m},1,r^u_{\ell,m})$. The
uniform expansion of $Df^{\pi(A)}$ and of $Df^{\pi(B)}$ in the
$E^{\uut}$ direction in the sets $U_A$ and $U_B$ imply that
$|r^u_{\ell,m}|\to 0$ as $\ell,m\to \infty$. Similarly, the
uniform contraction of $Df^{\pi(A)}$ and of $Df^{\pi(B)}$ in the
$E^{ss}$ direction gives $|r^s_{\ell,m}|\to 0$ as $\ell,m\to
\infty$.

Since the transitions $T^{\ct}_{(a,b)}$ and $T^{\ct}_{(b,a)}$ are
isometries and  central direction is preserved by $f_\nu$, the
central multiplier $\la^{\ct}(R_{\ell,m})$ of $R_{\ell,m}$
satisfies
$$
|\la^{\ct}(R_{\ell,m})|= |(\la^{\ct}(B))^m\, (\la^{\ct}(A))^\ell|
.
$$
Finally, by construction, the whole orbit of $R_{\ell,m}$ is
contained in the neighborhood $\neig$ of the simple cycle.
This completes the proof of the lemma.
{\hfill$\Box$\medskip}

Consider a partially hyperbolic saddle $P$ of a diffeomorphism $f$
with a splitting $E^{\sst}\oplus E^{\ct}\oplus E^{\uut}$, where
$E^{\ct}$ is a one-dimensional (central) direction. Consider the
strong stable manifold $W^{\sst}(P)$ of $P$ tangent to the strong
stable direction $E^{\sst}$ and the strong unstable manifold
$W^{\uut}(P)$ of $P$ tangent to $E^{\uut}$. Note that if $P$ is
hyperbolic and $E^{\ct}$ is expanding (resp. contracting) then
$W^{\sst}(P)=W^{\st}(P)$ (resp. $W^{\uut}(P)=W^{\ut}(P)$).

Given partially hyperbolic saddles $P$ and $Q$ as above and an
open set $U$ containing the orbits of $P$ and $Q$, for $\ri,\rj\in
\{\st,\sst,\ut,\uut\}$, we write $P\, \cap_{\ri,\rj,U}\, Q$ if
there is $X\in W^{\ri}(P)\cap W^\rj(Q)$ whose orbit is contained
in $U$ (i.e., $W^{\ri}(P)$ and $W^\rj(Q)$ are $U$-related). We
write $P\, \pitchfork_{\ri,\rj,U}\, Q$ if this intersection is
transverse.

\medskip

Note that the proof of Lemma~\ref{l.Rlm} immediately implies the
following:

\blm\label{r.I1}
The saddle $R_{\ell,m}$ in Lemma~\ref{l.Rlm} satisfies
$$
{\mbox{\bf (i)}}\quad A\, \pitchfork_{\st, \uut,\neig} \,
R_{\ell,m} \quad \mbox{and} \quad {\mbox{\bf (ii)}}\quad
R_{\ell,m}\, \pitchfork \, _{\sst,\ut,\neig} B.
$$
 \elm

\subsection{Simple cycles of biaccumulated saddles}
\label{ss.simplebiaccumulated}

In this section, we consider simple cycles associated to saddles
which are biaccumulated. We will see that, in this case, the
simple cycle can be chosen (after a
perturbation) satisfying some additional properties that we
proceed to explain.

Consider a  simple cycle associated to saddles $A$ and $B$, with
$\mbox{\rm index }(A)+1 = \mbox{\rm index }(B)$, such that $A$ is
$s$-biaccumulated. Let $X$ be the transverse heteroclinic point
of this cycle, $X\in W^{\st}(A)\cap W^{\ut}(B)$, and denote by $W^{\st}_X(A)$
the connected component of $(W^\st(A)\setminus W^{\sst}(A))$ containing
 $X$. By the $\st$-biaccumulation property, there is
some transverse homoclinic point $\zeta$ of $A$ in $W^{\st}_X(A)$.
 Then, by the Smale's homoclinic theorem, there is a  small
neighborhood $U_{A,\zeta}$ of the orbits of $A$ and $\zeta$ such
that the set $\La_f(U_{A,\zeta})$ is hyperbolic. Moreover, since
the saddle $A$ has real multipliers, after a perturbation of $f$,
we can assume that $\La_f(U_{A,\zeta})$ has a hyperbolic splitting
consisting  of one-dimensional bundles.
 In this case, we say that $U_{A,\zeta}$ is an {\emph{$\st$-adapted
 neighborhood of $A$ and $\zeta$}.

If the saddle $B$ is $\ut$-biaccumulated, there is
a transverse homoclinic point $\vartheta$ of $B$ in
the component $W^{\ut}_X(B)$ of  $(W^\ut(B)\setminus W^{\uut}(B))$ containing $X$. We define
{\emph{$\ut$-adapted
 neighborhoods of $B$ and $\vartheta$} in a similar way.

We use the next lemma to get relative cycles associated to $A$ and
$B$ in a set where the dynamics is partially
hyperbolic.

\blm
Consider a  simple cycle associated to saddles $A$ and $B$,
with $\mbox{\rm index }(A)+1 = \mbox{\rm index }(B)$, such that $A$
and $B$ are $\st$- and $\ut$-bi\-accu\-mu\-la\-ted, respectively. Let
$V_0$ be the neighborhood of the simple cycle and $X\in W^{\st}(A)\pitchfork
W^\ut(B)$ and $Y\in W^\ut(A)\cap W^\st(B)$ be its heteroclinic orbits.

Then there is $g$ arbitrarily close to $f$ with a simple cycle
associated to $A$ and $B$, heteroclinic points $X$ and $Y$, and an
associated neighborhood $V^\prime\subset V_0$, having the following
additional properties. There are transverse homoclinic points
$\zeta\in W^{\st}_X(A,g)$ of $A$ and $\vartheta\in W^{\ut}_X(B,g)$ of
$B$, disks $\Delta^\ut(\zeta)\subset W^\ut(A,g)$ and
$\Delta^\st(\vartheta)\subset W^\st(B,g)$, and adapted neighborhoods
$U_{A,\zeta}$ and $U_{B,\vartheta}$ such that
\begin{itemize}
\item
the maximal invariant set in
\begin{equation}
\label{e.V0}
 V= U_{A,\zeta}\cup U_{B,\vartheta}\cup V^\prime
\end{equation}
 admits a
$Dg$-invariant splitting consisting of one dimensional bundles,
\item in the  local coordinates in $A$ and after replacing $X^+$
 and $\zeta$ by forward
iterates, we have that $X^+=(0^s,-1,0^u)$,
$\zeta=(\zeta^s,-1,0^u)$, and
\begin{equation}\label{e.du}
\Delta^\ut(\zeta)=\{(\zeta^s,-1)\} \times [-1,1]^u \subset W^\ut(A,g)
\quad \mbox{and}\quad \bigcup_{i=0}^\infty g^{-i}(\Delta^\ut)\subset
V,
\end{equation}
\item
in the  local coordinates in $B$ and after replacing $X^-$
 and $\vartheta$ by backward
iterates, we have that $X^-=(0^s,1,0^u)$, $\vartheta=(0^s,
1,\vartheta^u)$, and
\begin{equation}
\label{e.ds}
 \Delta^\st(\vartheta)=[-1,1]^s\times
\{(1,\vartheta^u)\} \subset W^\st(B,g) \quad \mbox{and}\quad
\bigcup_{i=0}^\infty g^{i}(\Delta^\st)\subset V.
\end{equation}
\end{itemize}
\label{l.adapted}
 \elm

Observe that we can assume that the neighborhood $V$ is a
neighborhood satisfying Lemma~\ref{l.partialhyperbolicity}.

In \cite{BD4} (see, for instance, its Section 5.2)
are obtained similar results in a slightly different context. Thus
we just sketch the standard proof of the lemma.

\begin{proof} Since the dynamics in the sets $U_{A,\zeta},
U_{B,\vartheta}$,  and $V_0$ is partially hyperbolic, the first
assertion (partial hyperbolicity) immediately follows.

The proof of the second item consists of two steps. First, one
considers a small disk $\Delta^\ut_0\subset W^\ut(A)\cap
V_0$
containing some forward iterate of $\zeta$. After replacing
$\zeta$ and $\Delta^\ut_0$ by some forward iterates of them and a
perturbation, one can assume that $\Delta^\ut_0$ is parallel to the
unstable direction, that is, in the local coordinates at $A$, one
has
$$
\zeta=(\zeta^s,\zeta^c,0^u), \, \zeta^s\ne 0^s,\quad \mbox{and} \quad
\Delta^\ut_0=\{ (\zeta^s,\zeta^c)\} \times [-1,1]^u.
$$
Replacing $\zeta$ and $X^+$ by some iterates, we can assume that
$X^+=(0^s,-1,0^u)$ and $\zeta^c\in [-1,-\la^\ct(A))$. If $\zeta^c=-1$
we just take  $\Delta_0^\ut=\Delta^\ut(\zeta)$.  Otherwise, we
select a small neighborhood of $\zeta$ (which does not intersect
the central direction) and along a segment of orbit of $\zeta$ (we
will precise its length) we consider a perturbation $f_t$ (small
$t$) of $f$ which is a scaled $t$-translation of $f$ in the central direction. More
precisely, for simplicity let $\la=\la^\ct(A)$, then  in a vicinity
of $f_t^i(\zeta)$ the perturbation is of the form
$$
f_t(\xi)=f(\xi)+(0^s,\la^i\, t,0^u).
$$
We claim that for every large $k$ there is $t_k$, $t_k\to 0$ as
$k\to \infty$, such that
\begin{equation}
\label{e.tk} (f^{k}_{t_k}(\zeta))^c=-\lambda^{k+1}=(f^{k+1}_{t_k}
(X^+))^c=(f^{k+1} (X^+))^c.
\end{equation}
This implies that $f_{t_k}^{k} (\zeta)$ and $f_{t_k}^{k+1}(X^+)$ have
the same central component. Moreover, as the perturbation
preserves the unstable direction, after replacing $\zeta$ and
$X^+$ by $f_{t_k}^{k}(\zeta)$ and $f_{t_k}^{k+1}(X^+)=
f^{k+1}(X^+)$, we obtain the result for $g=f_{t_k}$. In that case, we consider a
perturbation only along the first $k$ iterates of $\zeta$.

To get the claim we first observe that one inductively gets
$$
(f_t^n(\zeta))^c= \la^n\, \zeta^c+ n\, \la^{n-1}\,t.
$$
Therefore, to get (\ref{e.tk}) it is enough to choose
$$
\la^{k}\, \zeta^c+ k\, \la^{k-1}\,t_k=-\la^{k+1}, \quad
t_k=\frac{-\la\,\zeta^c-\lambda^2}{k}.
$$

The proof of the third item is analogous to the one of the second one.
This completes the sketch of the proof of the lemma.
\end{proof}

\subsection{Quotient dynamics and heteroclinic/ho\-mo\-cli\-nic
intersections}\label{ss.intersections}  In this section, we study
the semi-local dynamics of diffeomorphism with simple cycles with
biaccumulation properties satisfying Lemma~\ref{l.adapted} (by
semi-local we mean the dynamics in a partially hyperbolic
neighborhood $V$ as in (\ref{e.V0})). The first step is to write the
homoclinic/heteroclinic intersection properties in
Proposition~\ref{p.bdf} (items {\bf 3)}--{\bf 5)}) in terms of the quotient dynamics.

Next lemma  (corresponding to \cite[Proposition 3.8]{BD4}) states a
relation between $V$-related intersections of  invariant manifolds
and the quotient dynamics. Write $\la=\la^\ct(A)$ and $\be=\la^\ct(B)$
and denote by $f_\lambda$ and $f_\beta$ the restrictions of $f^{\pi(A)}$ and $f^{\pi(B)}$
 to
the central directions in the neighborhoods $U_A$ and $U_B$ of the
saddles in the cycle. Note that these maps are just multiplications by $\la$ and
$\be$. In this notation, the subscript $\lambda$ and $\beta$ denote
the eigenvalue of this linear map.

\blm \label{l.Vintersections} Let  $f$ be a diffeomorphism with a
simple cycle with biaccumulation properties as in
Lemma~\ref{l.adapted}. Consider the parameter $\nu=\nu_{\ell,m}$,
the diffeomorphism $f_{\nu_{\ell,m}}$, and the saddle $R_{\ell,m}$
in Lemma~\ref{l.Rlm}. Assume that there are large $k$ and $h\in
\NN$,  $(k,h)\ne (\ell,m)$, such that
$$
f_\beta^h \circ T^\ct_{(a,b),\nu_{\ell,m}} \circ f_\la^k \circ
T^\ct_{(b,a)}(1)= f_\be^h(\tau\, f_\la^r(-1)+\nu_{\ell,m})= \be^h\,(
-\tau\, \la^r +\nu_{\ell,m})= 1.
$$
Then, if $\neig$ is the neighborhood in \eqref{e.V0}, the
following holds
$$
\begin{array}{lll}
&{\mbox{\bf (i)}}\quad
R_{\ell,m} \, \cap_{\sst,\ut,\neig}\,
 A,
&\quad
{\mbox{\bf (ii)}}\quad
B \, \cap_{\st,\uut,\neig} R_{\ell,m},
\\
&{\mbox{\bf (iii)}}\quad
R_{\ell,m} \, \cap_{\sst,\uut,\neig}\,R_{\ell,m},
&\quad{\mbox{\bf (iv)}}\quad
B \, \cap_{\st,\ut,\neig} \, A.
\end{array}
$$
 \elm

\brm \label{r.Vintersections}
The intersection properties in Lemmas~\ref{l.Vintersections} and \ref{r.I1}
imply items 3)--5) in
Proposition~\ref{p.bdf}:
\begin{itemize}
\item
item 3) ($B$ and $R_{\ell,m}$ are homoclinically related)
follows immediately from ii) in Lemma~\ref{r.I1} and ii) in Lemma~\ref{l.Vintersections};
\item
the $V$-related intersections in item 4) are just iii) and ii) in Lemma~\ref{l.Vintersections};
\item
item 5) (cycle associated to $A$ and $R_{\ell,m}$) follows from
i) in Lemma~\ref{r.I1} and i) in Lemma~\ref{l.Vintersections}.
\end{itemize}
\erm

\begin{proof}
The proof follows arguing exactly   as in the construction of the points $R_{\ell,m}$
in Lemma~\ref{l.Rlm}. We explain how to obtain the relation
$R_{\ell,m}\, \cap_{\uut,\sst,V_0} \, R_{\ell,m}$. The other properties follow
analogously after doing the corresponding identifications via the
quotient by the strong stable/unstable hyperplanes. In the
coordinates in $U_B$ we have
$$
W^{\uut}_{\loc}(R_{\ell,m})=\{(r^s_{\ell,m},1)\} \times [-1,1]^u.
$$
Since the directions $E^{\sst}$, $E^\ct$ and $E^{\uut}$ are preserved, by
the definition of $(k,h)$,
there is a disk $\De^\ut\subset
[-1,1]^u$ such that
$$
\begin{array}{ll}
& f^{h\, \pi(B)+t_{(a,b)}+ k\, \pi(A) + t_{(b,a)}}_{\nu_{\ell,m}} (
\{(r^s_{\ell,m},1)\} \times \De^\ut)= \\
&\qquad \quad= f^{h\,\pi(B)}\circ \mfT_{(a,b),{\nu_{\ell,m}} } \circ f^{k\,\pi(A)}
\circ
\mfT_{(b,a)} ( \{(r^s_{\ell,m},1)\}\times \De^\ut)=\\
&\qquad\quad =\{ (\xi^s,f_\be^h\circ T^\ct_{(a,b),\nu_{\ell,m}} \circ
f_\la^k \circ
T^\ct_{(b,a)}(1))\} \times [-1,1]^u=\\
 &\qquad\quad=\{(\xi^s,1)\} \times [-1,1]^u,
\end{array}
$$
for some $\xi^s$. Since, in the coordinates in  $U_B$, $ W^{\sst}_{\loc}
(R_{\ell,m})=[-1,1]^s\times \{(1,r_{\ell,m}^u)\}$, we get that
$W^{\uut}(R_{\ell,m})$ meets $W^{\sst}(R_{\ell,m})$. Finally,  this
intersection can be taken $V_0$-relative: we just consider points
of the disk $\{(r^s_{\ell,m},1)\} \times \De^\ut$ whose (segment) of
orbit remains in $\neig$. Actually, in this part we do not use the
biaccumaltion properties and only consider iterations in a
neighborhood of the simple cycle.

To get the other properties in the lemma, we proceed exactly as
above considering the disks $\Delta^\ut(\zeta)=\{(\zeta^s,-1)\}
\times [-1,1]^u\subset W^\ut(A)$ and
$\Delta^\st(\vartheta)=[-1,1]^s\times \{(1,\vartheta^u)\}\subset
W^\st(B)$ in Lemma~\ref{l.adapted} (see equations \eqref{e.du} and
\eqref{e.ds}).
 This completes the sketch of the proof of
Lemma~\ref{l.Vintersections}.
\end{proof}

\medskip

We now see how the conditions in Lemma~\ref{l.Vintersections} are
obtained in our setting. For that we modify the maps
$f_\la$ and $f_\be$
 locally, but the resulting diffeomorphisms still
preserve the bundles $E^{\sst}, E^\ct$ and $E^{\uut}$.

\blm {\bf{(\cite[Corollaries 3.13 and 3.15]{BD4})}} We use the
notation above. For every large $K\in \NN$ and every $\ve>0$ there
are $\bar\beta \in (\beta-\ve,\beta+\ve)$, $\xi\in (0,\ve)$, and
natural numbers $k,p,q>K$, $k$ is even,  such that
$$
f_{\bar \beta}^p \circ T^\ct_{(a,b),\nu_k}\circ f_{\la}^k \circ
T^\ct_{(b,a)}(1)=1 \quad \mbox{and} \quad f_{\bar \beta}^q \circ
T^\ct_{(a,b),\nu_k}\circ f_{\la}^{k-2} \circ T^\ct_{(b,a)}(1)=1,
$$
where
$$
\nu_k=\left\{
\begin{array}{lll} &|\la^{k-2}|+\xi, &\quad \mbox{if $T^{\ct}_{(a,b)}$
preserves the orientation,}\\
 & -|\la^{k}|+\xi, &\quad \mbox{if $T^{\ct}_{(a,b)}$
reverses the orientation}.
\end{array}\right.
$$
Moreover, if $T^\ct_{(a,b)}$ preserves the orientation then $p=p(k)$ and $q=q(\xi)$,
$q(\xi)\to \infty$ as $\xi\to 0^+$,
and if $T^\ct_{(a,b)}$ reverses the orientation then $q=q(k)$ and $p=p(\xi)$
($p(\xi)\to \infty$ as $\xi\to 0^+$).
\label{l.corbd}
\elm
\begin{proof}
To give an idea of the proof of Lemma~\ref{l.corbd}
we consider the  orientation preserving case
(the orientation reversing case is similar and is left to a reader).
Note that if we take $\hat \nu=|\la^{k-2}|$,
for some large $k$, then there are $p=p(k)$ and $\hat \beta$ close
to $\beta$ with
$$
\hat\beta^{p}\, (-|\la^k| +|\la|^{k-2})=1.
$$
Therefore
$$
f_{\hat \beta}^p \circ T^\ct_{(a,b),\hat \nu}\circ f_{\la}^k \circ
T^\ct_{(b,a)}(1)=\hat\beta^{p}\, (-|\la^k| +|\la|^{k-2})=
1.
$$
Note also that this choice of $\hat\nu$ gives
$$
T^\ct_{(a,b),\hat \nu}\circ f_{\la}^{k-2} \circ
T^\ct_{(b,a)}(1)=(-|\la^{k-2}| +|\la|^{k-2})=
0.
$$

\begin{figure}[htb]

\begin{center}
\psfrag{bp}{$\hat{b}^{-p}$}
\psfrag{lk}{$\la^k$}
\psfrag{lk2}{$\la^{k-2}$}
\psfrag{lkn}{$\la^{k}+\nu$}
\psfrag{n}{$\hat\nu$}
\psfrag{x}{$\hat\nu+\xi$}
\psfrag{s}{$\xi$}
\psfrag{0}{$0$}
\psfrag{1}{$1$}
\psfrag{bq}{$\bar{b}^{-q}$}
\psfrag{dp}{$\bar{b}^{-p}$}
\includegraphics[height=2.2in]{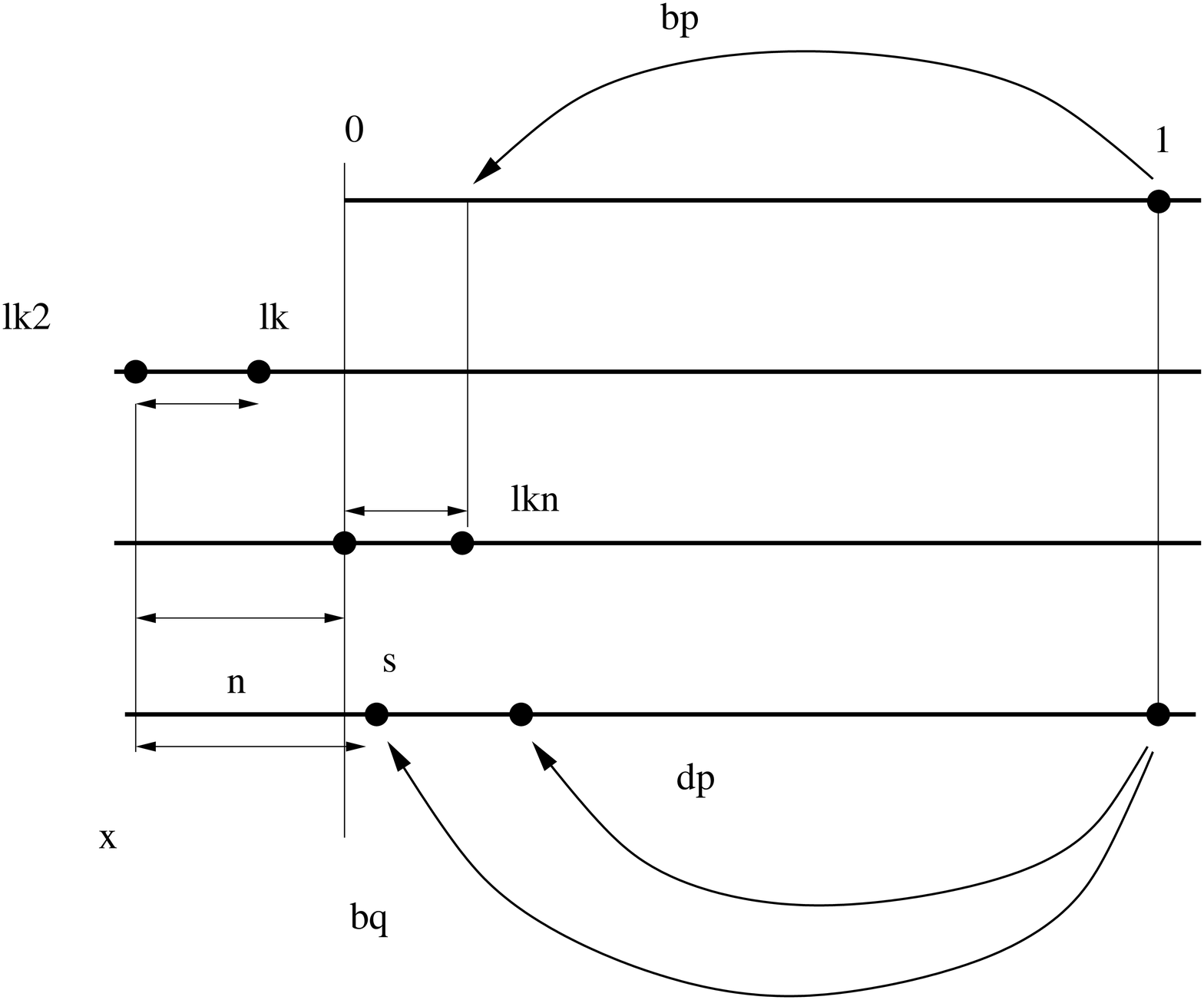}

 \caption{}
\label{f.fixedpoints}
\end{center}

\end{figure}

To complete the proof we just need to choose arbitrarily small $\xi$
and a new $\bar\beta$ (close to $\hat\beta$) such that
$$
\bar\beta^{-q}=\xi
\quad
\mbox{and}
\quad
\bar\beta^{p}\, (-|\la^k| +|\la|^{k-2}+\xi)=1.
$$
Then for
$\nu=\hat\nu+\xi$ the following equations hold simultaneously (see
Figure~\ref{f.fixedpoints}):
$$
\bar\beta^{-q}=\xi \quad  \Rightarrow  \quad
f_{\bar \beta}^q\circ
T^\ct_{(a,b),\hat\nu+\xi}\circ f_{\la}^{k-2} \circ
T^\ct_{(b,a)}(1)=1
$$
for some $q=q(\xi)$, $q(\xi)\to \infty$ as $\xi\to 0$,
and
$$
\bar\beta^{p}\, (-|\la^k| +|\la|^{k-2}+\xi)=1
\quad  \Rightarrow  \quad
f_{\bar \beta}^p \circ T^\ct_{(a,b),\hat \nu+\xi}\circ f_{\la}^k \circ
T^\ct_{(b,a)}(1)=1.
$$
This completes the proof of the lemma.
\end{proof}

Notice that, using the notation above, there are two periodic
points for $f_\nu$, $R_{k,p}$ and $R_{k-2,q}$.

We now explain our
choice of the periodic point $R_{\ell,m}$. This choice depends on whether or not
the map $T^\ct_{(a,b)}$ preserves the orientation, and
is done bearing in mind that we want to  get periodic points $R_{\ell,m}$ with uniformly
bounded (with respect to  $\ell$ and $m$ chosen in a specified way) central multipliers.

 First,
if $T^{\ct}_{(a,b)}$ is orientation preserving then we let $\ell=k$
and $m=p$. To calculate the central multiplier of $R_{\ell,m}$ note that
$\la^\ell=|\la^\ell|$ ($\ell=k$ is even) and that, in this case,
 $f_{\bar \beta}^m \circ
T^\ct_{(a,b),\nu}\circ f_{\la}^\ell \circ T^\ct_{(b,a)}(1)=1$ for $\nu=
|\la^{\ell-2}|+\xi$ means

\begin{equation}\label{e.Rlm+}
\bar\beta^{m}\, (-|\la^\ell| +|\la|^{\ell-2}+\xi)=1 \quad
\Rightarrow \quad |\bar \be^m\,\la^\ell|=\frac{(1-\bar
\be^m\,\xi)\, \la^2}{1-\la^2}. \end{equation}
In the reversing
orientation case, we let $\ell=k-2$ and $m=q$. In this case,
$f_{\bar \beta}^m \circ T^\ct_{(a,b),\nu}\circ f_{\la}^\ell \circ
T^\ct_{(b,a)}(1)=1$ for $\nu=-|\la^{\ell+2}|+\xi$ implies that
\begin{equation}\label{e.Rlm-}
\bar\beta^{m}\, (|\la^\ell| -|\la|^{\ell+2}+\xi)=1 \quad
\Rightarrow \quad
|\bar \be^m\,\la^\ell|=\frac{(1-\bar\be^m\,\xi)}{1-\la^2}.
\end{equation}

Notice that here we used essentially the specific choice of $\ell$ and
$m$ above. Notice also that in both cases the choice of $\ell$ and
$m$ can be done before the choice of $\xi$.
Thus, for chosen $m$ and $\ell$, taking small $\xi$ implies that the absolute value of the central multiplier of $R_{\ell,m}$  (which is $|\bar \be^m\,\la^\ell|$)
is uniformly bounded (independently of chosen $\ell$ and $m$). Finally, we get the following statement.

\blm
\label{r.uniform}
There is $\const>0$ such that
the absolute value of the  central multiplier of the point $R_{\ell,m}$ satisfies
${\const}^{-1}< |\la^\ct(R_{\ell,m})|<\const$,
for all large $\ell$ and $m$ chosen as above.
Therefore, the saddle $R_{\ell,m}$ satisfies item 2) in Proposition~\ref{p.bdf}.
\elm

\subsection{Proofs of Propositions~\ref{p.bdf} and
\ref{p.simpleanddense}} \label{ss.proofs}

\subsubsection{Proof of Proposition~\ref{p.bdf}}\label{sss.pbdf}
Consider the partially hyperbolic neighborhood $V$ in (\ref{e.V0})
in Lemma~\ref{l.adapted} and the sequence of parameters $\nu_k$ in
Lemma~\ref{l.corbd}. Let $R_k$ be the periodic point of $f_{\nu_k}$
given by Lemma~\ref{l.Vintersections} (with the notation of this
lemma, $R_k=R_{\ell_k,m_k}$). We claim that these saddles satisfy
the conclusions in Proposition~\ref{p.bdf}. The period of $R_k$ is
given by Lemma~\ref{l.Rlm} and the estimate on the central
multiplier and the itinerary of this saddle are
 provided by Lemma~\ref{r.uniform}. Finally, since the
splitting of $\La_{f_k}(V)$ consists of one-dimensional directions
(Lemma~\ref{l.adapted}) this saddle has real multipliers.

Taking small $\xi>0$ in
 equations
\eqref{e.Rlm+} and \eqref{e.Rlm-} and noting that $m$ remains fixed
and $\la=\la^\ct(A_f)$ is close to one, one gets
that
$R_k$ and $B_f$  have the same index.

The intersection properties in items 3)--5) follow from Lemmas~\ref{r.I1} and \ref{l.Vintersections},
as stated  in Remark~\ref{r.Vintersections}.
 This concludes the
proof of Proposition~\ref{p.bdf}.
\hfill $\Box$

\subsubsection{Proof of Proposition~\ref{p.simpleanddense}}
\label{sss.psimpleanddense}

The proof of Proposition~\ref{p.simpleanddense} follows from the
previous constructions and the ones in \cite{BD4}, where it is
proved that every co-index one cycle generates robust cycles
\footnote{Consider  $f$ with  a co-index one cycle, then there is
a $C^1$-open set $\cU$ whose closure contains $f$ such that very
$g\in \cU$ has transitive hyperbolic sets $\La_g$ and $\Si_g$ of different
indices such that $W^u(\Si_g)\cap W^s(\La_g)\ne \emptyset$ and
$W^s(\Si_g)\cap W^u(\La_g)\ne \emptyset$. Necessarily, at least
one of these sets is non-trivial.}.
The key step of the
construction in \cite{BD4} is to obtain a partially hyperbolic saddle $S$ with a real
multiplier close to one having a strong stable/unstable connexion
$W^{\sst}(S)\cap W^{\uut}(S)\ne\emptyset$, exactly as the points
$R_{k}$ in Proposition~\ref{p.bdf} (note that,
{\emph{a priori,}} the central multiplier of $R_k$ is not close to one).
Such a dynamical configuration generates  cycles via a
blender-like construction that we will explain below. As the existence
of a blender structure is open (see \cite[Lemma 1.11]{BD1}), we can explain the proof for a
special linear perturbation. We will follow the
model in \cite{BDV1}.
For a discussion of the notion of blender see \cite[Chapter 6.2]{BDV2}.
We now go to the details of this construction.

Applying Proposition~\ref{p.bdf} to the partially hyperbolic neighborhood $V$,
we get a sequence of diffeomorphisms $f_k$ converging to $f$ such that
every $f_k$ has
a saddle $R_k$ satisfying the conclusions of Proposition~\ref{p.bdf}. To get these
saddles having central multipliers  with modulus close  to one we use the
following lemma:

\blm[Franks' Lemma, \cite{F}]
\label{l.franks}
 Consider a diffeomorphism  $f$
 and an $\ve$-perturbation $A$ of the
derivative $Df$ of $f$ along an $f$-in\-va\-riant finite set
$\Sigma$. Then, for every neighborhood $U$ of $\,\Sigma$, there is a
diffeomorphism $g$ $C^1$-$\ve$-close to $f$ such that $g(x)=f(x)$,
if  $x\in \Sigma$ or if $x\not\in  U$, and $Dg(x)=A(x)$, for all
$x\in \Sigma$. \elm

Recall that the periods $\pi(R_k)$ go to infinity as $k\to \infty$ and that
the central multipliers $|\la^\ct(R_k)|$ are uniformly bounded. To get a saddle $R_k$
(satisfying Proposition~\ref{p.bdf})
whose
central multiplier has modulus close to $1$  it is enough to consider a multiplication
in the central direction along the orbit of $R_k$ by a factor $(1+\varepsilon)\,
|\la^\ct(R_k)|^{-1/\pi(R_k)}$
(which is close to $(1+\varepsilon)$ for large $k$) preserving the partially hyperbolic splitting
and apply Lemma~\ref{l.franks}.

\medskip

We now proceed to explain the blender-like construction when the
multiplier of $R_k$ is positive (see \cite{BD4} for the negative
case)\footnote{In fact, our construction can be modified to get the
saddles $R_k$ in Proposition~\ref{p.bdf} having positive central
multipliers, but this construction is much more involved.}. Assume
that, in our local coordinates, the point $R=R_{k}=(0^s,0,0^u)$ and that the local dynamics at
$R$ is linear, $f(x^s,x,x^u)=(A^\st(x^s), \sigma\, x, A^\ut(x^u))$,
where $A^\st$ is a contracting matrix, $A^\ut$ is expanding, and
$\la^\ct(R)=\sigma\in (1,2)$. By 4) in Proposition~\ref{p.bdf},
there is a strong stable/unstable intersection associated to $R$:
there are a small disk $D^\ut\subset \{(0^s,0)\} \times
[-1,1]^u\subset W^\ut(R)$ and $k_0$ such that $f^{k_0}(D^\ut)=
\{(\kappa^s_0, 0)\} \times [-1,1]^u$ (the first $k_0$ iterates of
$D^\ut$ are contained in $\neig$). We  consider a one parameter family
of diffeomorphisms (perturbations of $f$ in a neighborhood of $D^\ut$)
such that if $(x^s,x,x^u)$ is in a neighborhood of $D^\ut$ then
$$
f_t^{k_0}(x^s,x,x^u)=f(x^s,x,x^u)+(0^s,t,0^u).
$$
Consider $t<0$ and the cube
$$
C=[-1,1]^s\times [0,|t|/(\sigma-1)]\times [-1,1]^u\subset V
$$
 Now \cite[lemma in page 717]{BDV1} implies the following intersection
 property:

\medskip

\noindent {\bf Intersection property:}
{\emph{Consider the disk  $\De^\ut_{[\rho_1,\rho_2]}=\{d_0^s\} \times [\rho_1,\rho_2]\times [-1,1]^u\subset C$,
$\rho_1<\rho_2$.
 Then there is a point $X$ of transverse
intersection between
 $W^{\st}(R)$ and $\De^\ut_{[\rho_1,\rho_2]}$ whose forward orbit is contained in $V$,
see Figure~\ref{f.intersection}.}

\begin{figure}[htb]

\begin{center}
\psfrag{R}{$R$} \psfrag{D}{$\De^\ut_{[\rho_1,\rho_2]}$ }
\psfrag{td}{$|t|/(\sigma-1)$} \psfrag{Ws}{$W^\st(R)$} \psfrag{t}{$t$}
\includegraphics[height=1.4in]{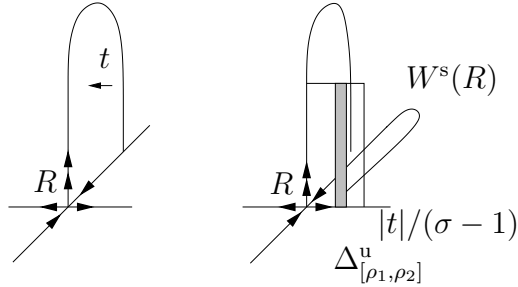}
 \caption{Intersection properties}
\label{f.intersection}
\end{center}

\end{figure}

\medskip


Since we  have $B\, \cap_{\st,\uut,V}\, R$ and $A\, \cap_{\ut,\st,V}\, R$
(recall items 4) and 5) of  Proposition~\ref{p.bdf}), there are
compact disks $K^\st(B)\subset W^\st(B)$ and $K^\ut(A)\subset W^\ut(A)$
(contained in $V$) close to
 $W^{\st}_{\loc}(R)$
and  $W^{\uut}_{\loc}(R)$, respectively.
We can now consider a local perturbation $g$ of $f_t$  (in a domain different of the family
of perturbations $f_t$)
such that in the local coordinates at $R$ one has (see Figure~\ref{f.kuks})
$$
\begin{array}{ll}
K^\ut(A)&=\{(a_0^s, a_0)\}\times [-1,1]^u, \quad a_0\in
(0,|t|/(\sigma-1));\\
K^\st(B)&=[-1,1]^s\times \{ (b_0, b^u_0)\}, \quad b_0<0.
\end{array}
$$


\begin{figure}[htb]

\begin{center}
\psfrag{R}{$R$}
\psfrag{ku}{$K^u(A)$ }
\psfrag{ks}{$K^s(B)$}
\psfrag{t}{$t$}
\includegraphics[height=1.4in]{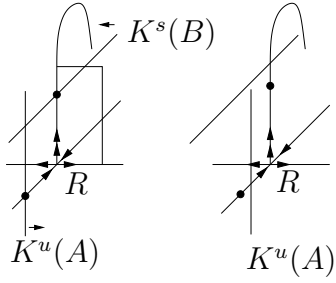}
 \caption{The perturbation $g$}
\label{f.kuks}
\end{center}

\end{figure}

\medskip

Finally, the disk $K^\ut(A)$ intersects the cube $C$ from the bottom to the top.
And this holds for the continuation of $K^\ut(A_h)$ for every $h$ close to $g$.
This implies that there is a  $C^1$-neighborhood $\cE$ of $g$ such that
(by the intersection property) $K^\ut(A_h)$ is $V$-accumulated by the stable manifold
of $R_h$. Therefore densely in $\cE$ one can obtain
$V$-related cycles associated
to $R_h$ and $A_h$. Since $B_h$ is $V$-homoclinically related to $R_h$,
one also has $V$-related cycles associated to $B_h$ and $A_h$. This completes the
proof of Proposition~\ref{p.simpleanddense}. \hfill $\Box$

\section{Final construction}  \label{s.genzero}



 In this section we use the results of previous sections to complete the proof of Theorem~\ref{t.main}.

Let $\mathcal{G}$ be the residual set of $\text{\rm Diff}^{\,1}(M)$
described in Section~\ref{ss.generic}. The statement below is a local version of Theorem~\ref{t.main}.
\bthm\label{t.two} Let $f\in \mathcal{G}$  have
two hyperbolic saddles $P_f$ and $Q_f$ of different indices in the same homoclinic class.  Arbitrarily $C^1$-close to
$f$ there exists a $C^1$-open set $\mathcal{Z}\subset \text{\rm
Diff}^{\,1}(M)$ and a residual subset $\mathcal{R}\subset
\mathcal{Z}$ such that every $g\in \mathcal{R}$ has a non-hyperbolic
ergodic invariant measure with uncountable support inside of the
homoclinic class of  $P_g$ and $Q_g$.
\ethm
\brm
The genericity  hypothesis imply that $H(P_g, g)=H(Q_g, g)$.
\erm
Due to
standard genericity arguments Theorem~\ref{t.two} implies Theorem~\ref{t.main},
for details see \cite{ABCDW}.

\medskip

The genericity hypothesis imply that the map $f$ 
has two saddles with real multipliers of consecutive indices,
say $A_f$ and $B_f$, $\text{\rm index}(A_f)+1=\text{\rm
index}(B_f)$, such that the saddle $A_f$ is homoclinically related to $P_f$ and the saddle
$B_f$ is homoclinically related to $Q_f$. Moreover, we can assume that the saddle $A_f$ is
$\st$-biaccumulated, the saddle $B_f$ is $\ut$-biaccumulated,
and that the central multiplier of $A_f$ is close to one (it is enough to have $|\lambda^{\ct}(A_f)|\in (0.9, 1)$),
see generic conditions {\bf R3)}--{\bf R4)}.

We proceed as follows:

\begin{itemize}
 \item
Proposition~\ref{p.bdhayashi} implies that  by a
$C^1$-perturbation we can create a cycle
corresponding to the continuations of the saddles $A_f$ and $B_f$. 
\item
 Apply Lemma~\ref{l.partialhyperbolicity} to this cycle. This gives
 a map $g$ which is $C^1$-close to $f$, and an open set $V\subseteq
 M$ such that $\Lambda_g(V)$ is strongly partially hyperbolic
(one dimensional central manifold), $A_g$ is ${V}$-$\st$-biaccumulated and $B_g$ is
${V}$-$\ut$-biaccumulated, and $g$ has a $V$-related cycle associated
to $A_g$ and $B_g$.
\item
By Proposition~\ref{p.simpleanddense}, $C^1$-near
there is an open set $\mathcal{Z}$ and a dense in $\mathcal{Z}$
countable subset $\mathcal{D}\subset \mathcal{Z}$ such that every
$g\in \mathcal{D}$ has a  $V$-related cycle associated with
the saddles $A_g$ and $B_g$.
\end{itemize}

 From now on we fix the open set $\mathcal{Z}\subset \text{\rm Diff}^{\,1}(M)$,
the countable dense subset $\mathcal{D}$ of $\mathcal{Z}$, the saddles $A_g$ and $B_g$,
 and the neighborhood $V\subset M$.

 \bpro\label{p.localperiod}
  Generic diffeomorphisms from $\mathcal{Z}$ have a
sequence of periodic saddles in $V$  which satisfies the assumptions of
Proposition~\ref{p.firstmainproposition} and belongs to the continuation of the homoclinic class above.
 \epro

Note that Propositions~\ref{p.localperiod} and \ref{p.firstmainproposition} imply
the existence of non-hyperbolic ergodic measures for generic diffeomorphisms from $\mathcal{Z}$, and, thus, Theorem~\ref{t.two}.

\medskip

 For any $g\in \mathcal{Z}$ the maximal invariant set
$\Lambda_g=\cap_{k\in \mathbb{Z}}g^k(V)$  is a (not necessarily closed) partially hyperbolic invariant set (with a splitting
$E^{\sst}\oplus E^\ct \oplus E^{\uut}$, $E^\ct$ is one-dimensional).
 For $g\in \mathcal{Z}$ denote by $\chi^{\ct}(A_g)<0$
and $\chi^{\ct}(B_g)>0$ the central Lyapunov exponents of saddles
$A_g$ and $B_g$ (corresponding to the central direction $E^{\ct}$). Fix a constant $C$ such that
\begin{equation}
\label{e.C}
C>\sup_{g\in \mathcal{Z}}\displaystyle{\frac{16}{
|\chi^{\ct}(A_g)|}}.
\end{equation}

\bpro\label{p.conditions} For each $N\in \mathbb{N}$ there is a
family of open sets $ \mathcal{Z}_{n_1, \ldots ,n_N}\subseteq
\mathcal{Z} $ indexed by $N$-tuples $(n_1, \ldots ,n_N), \ n_i\in
\mathbb{N},$ satisfying the following properties:

\begin{description}

\item{{\bf Z1)}}
For  any tuples $(n_1, \ldots ,n_N)\ne (m_1, \ldots ,m_N)$ the sets
$ \mathcal{Z}_{n_1, \ldots ,n_N}$ and $\mathcal{Z}_{m_1, \ldots,
m_N}$ are disjoint.
\item{{\bf Z2)}}
For  any tuple $(n_1, \ldots, n_N, n_{N+1})$ we have $
\mathcal{Z}_{n_1, \ldots ,n_N, n_{N+1}}\subseteq \mathcal{Z}_{n_1,
\ldots ,n_N}. $ In particular, $\mathcal{Z}_{n_1}\subseteq
\mathcal{Z}$ for every $n_1\in \mathbb{N}$.
\item{{\bf Z3)}}
The union
$\cup_{n_1\in \NN}\mathcal{Z}_{n_1}$ is dense in $\mathcal{Z}$, and
for each $N$
the union $\cup_{j\in \mathbb{N}} \mathcal{Z}_{n_1, \ldots n_{N},j}$
is dense in $\mathcal{Z}_{n_1, \ldots ,n_N}$.
\item{{\bf Z4)}}
Every diffeomorphism $g\in \mathcal{Z}_{n_1, \ldots ,n_N}$ has a
finite sequence of periodic saddles homoclinically related to $B_g$
(thus of the same index as $B_g$)
$$
\{P_{n_1}, P_{n_1,n_2},\ldots , P_{n_1,n_2,\ldots ,n_N}\}\subset
\Lambda_g \cap H(B_g,g)
$$
having
 real multipliers, satisfying the  $V$-$\ut$-biaccumulation property,
and of growing periods, $\pi(B_g)<\pi (P_{n_1})<\ldots
<\pi(P_{n_1,n_2,\dots ,n_N})$. Moreover, saddles $\{P_{n_1}, P_{n_1,n_2},\ldots , P_{n_1,n_2,\ldots ,n_N}\}$ depend continuously on $g$ when $g$ varies over $\mathcal{Z}_{n_1, \ldots ,n_N}$. 
\item{{\bf Z5)}}
For  any tuple $(n_1, \ldots ,n_{N})$ there exists a countable dense
subset  $\mathcal{D}_{n_1, \ldots ,n_{N}}\subset \mathcal{Z}_{n_1,
\ldots ,n_{N}}$ such that  every $g\in \mathcal{D}_{n_1, \ldots,
n_{N}}$ has a $V$-related heterodimensional cycle associated to
the saddles $A_g$ and $P_{n_1,\ldots ,n_N}$. 
\item{{\bf Z6)}}
There are numbers $\{\gamma_{n_1, \ldots ,n_{N}}\}_{(n_1 ,\ldots,
n_{N})\in \mathbb{N}^N}$ such that for any $N\in \mathbb{N}$, any
tuple $(n_1 ,\ldots ,n_N, n_{N+1})$, and any $g\in \mathcal{Z}_{n_1,
\ldots ,n_N ,n_{N+1}}$ the orbit of $P_{n_1, \ldots,n_N, n_{N+1}}$
is a $(\gamma_{n_1, \ldots ,n_{N}}, 1- C\, |\chi^{c}(P_{n_1, \ldots,
n_N})|)$-good approximation  of the orbit of $P_{n_1, \ldots, n_N}$
(recall Definition~\ref{d.goodaprox}), where $C$ is the constant in
\eqref{e.C}.
\item{{\bf Z7)}}
Take any $g\in \mathcal{Z}_{n_1, \ldots ,n_{N}}$. Let $d_k, 1\le
k\le N$,  be the minimal distance between the points of the $g$-orbit of
$P_{n_1, \ldots ,n_{k}}$. Then
$$
\gamma_{n_1 ,\ldots, n_{N}}< \frac{\min_{1\le k\le N} d_k}{3\cdot 2^N}. 
$$
\item{{\bf Z8)}}
For any $N\in \mathbb{N}$, any tuple $(n_1, \ldots, n_N, n_{N+1})$,
and any $g\in \mathcal{Z}_{n_1, \ldots,n_N ,n_{N+1}}$
$$
|\chi^{c}(P_{n_1, \ldots, n_N, n_{N+1}})|< \frac{1}{2}\,
|\chi^{c}(P_{n_1, \ldots, n_{N}})|.
$$
\end{description}
\epro

Before proving Proposition~\ref{p.conditions} let us complete the
proof of Proposition~\ref{p.localperiod} (and, therefore, of the
main result).

Due to Property {\bf Z3)}, for any $N\in
\mathbb{N}$ the set $\widetilde{\mathcal{Z}}_{N}=\cup_{(n_1, \ldots
,n_{N})\in \mathbb{N}^N} \mathcal{Z}_{n_1, \ldots ,n_{N}}$ is an open
and dense subset of $\mathcal{Z}$. Consider the intersection
$\mathcal{R}=\cap_{N\in \mathbb{N}}\widetilde{\mathcal{Z}}_N$. The
set $\mathcal{R}$ is a residual subset of $\mathcal{Z}$. Take any
$g\in \mathcal{R}$. Property {\bf Z1)} implies that for each
$N\in \NN$ the map $g$ belongs to one and only one set from the
collection $\{\mathcal{Z}_{n_1, \ldots, n_N}\}_{(n_1, \ldots,
n_N)\in \NN^N}$. Therefore, due to {\bf Z2)} and {\bf Z4)}, for the
diffeomorphism $g\in \mathcal{R}$ a sequence of periodic points
$\{B_g, P_{n_1}, P_{n_1,n_2}\ldots , P_{n_1,n_2,\ldots, n_N}, \ldots
\}\subset \Lambda_g \cap H(B_g,g)$ is well defined.  We claim that this sequence
satisfies the assumptions of
Proposition~\ref{p.firstmainproposition}.

Indeed, assumptions {\bf
1)} and {\bf 2)} follows from the choice of the set $V\subset M$ and
Property {\bf Z4)}. Assumptions {\bf 3)}, {\bf 4)}, and {\bf 5)}
follow from {\bf Z6)},   {\bf Z7)},  and {\bf Z8)}, respectively.
This proves Proposition~\ref{p.localperiod}, and an application of
Proposition~\ref{p.firstmainproposition} now implies
Theorem~\ref{t.two} (and, hence, Theorem~\ref{t.main}).

\begin{proof}[Proof of Proposition~\ref{p.conditions}]
Let us first construct sets $\mathcal{Z}_{n_1}\subset \mathcal{Z},
n_1\in \mathbb{N}$. The subset $\mathcal{D}\subset \mathcal{Z}$
(consisting of diffeomorphisms with cycles)
is countable. Let us enumerate diffeomorphisms from  $\mathcal{D}=\{g_i\}_{i\in \mathbb{N}}$. Take one of these diffeomorphisms, say, $g_i\in \mathcal{D}$.

Denote by $\cO(P)$ the orbit of a point $P$.
Choose a constant $\gamma>0$ such that 

\begin{equation}
\label{e.gamma}
\gamma<\inf_{g\in
\mathcal{Z}} \ \ \min_{X, Y\in \mathcal{O}(B_g),\ X\ne
Y}\frac{1}{3}\,\text{\rm dist} (X,Y).
\end{equation}
 Denote by $\lambda<1$ the
central multiplier of the saddle $A_{g_i}$, and
by $\beta>1$ the central multiplier of the saddle $B_{g_i}$. Choose
a neighborhood $\mathcal{U}(g_i)\subset \text{\rm Diff}^{\, 1}(M)$
and small (in particular, each component is of radius smaller than
$\gamma$) neighborhoods $U_{A_{g_i}}$,
and $U_{B_{g_i}}\subset M$ of the orbits of $A_{g_i}$ and $B_{g_i}$ such that
\begin{equation}\label{e.lambda}
   \forall
g\in\mathcal{U}(g_i) \ \forall x\in U_{A_{g_i}}\colon \quad  \lambda^2\le \|Df^{\pi(A_{g_{i}})}|_{E^c}(x)\|\le \lambda^{\frac{1}{2}}
 \end{equation}
and
\begin{equation}\label{e.betabis}
   \forall g\in\mathcal{U}(g_i) \ \forall x\in U_{B_{g_i}}  \colon
\quad \beta^{\frac{1}{2}}\le \|Df^{\pi(B_{g_{i}})}|_{E^c}(x)\|\le \beta^{2}.
 \end{equation}

Proposition~\ref{p.bdf} allows to obtain a sequence of diffeomorphisms $g_{ik}$,
$g_{ik}\to g_i$ as $k\to \infty$, such that
each diffeomorphism $g_{ik}$ has a periodic saddle $S_{ik}$ with real multipliers (denoted by $R_k$ in Proposition~\ref{p.bdf}),
having the following properties:
\begin{description}
\item{\bf S1)}
the saddle $S_{ik}$ is $V$-homoclinically related to $B_{g_{ik}}$, thus has the same index as $B_{g_i}$;
\item{\bf S2)}
for some constant $\const$ that does not depend on $k$ and for each $k\in \mathbb{N}$ we have $1<|\lambda^c(S_{ik})|<\const$;
\item{\bf S3)}
the map $g_{ik}$ has a $V$-related 
cycle associated to $S_{ik}$ and $A_{g_{ik}}$;
\item{\bf S4)}
there are  sequences of natural numbers $\ell_k,m_k$ that tend to infinity as $k\to \infty$,
such that under the iterates of $g_{ik}$ the saddle $S_{ik}$ needs
 a fixed number of iterates 
 (independent of $k$) to go from a  neighborhood $U_{B_{g_i}}$
 to a neighborhood $U_{A_{g_i}}$, then it remains $\ell_k\,\pi(A_{g_i})$ iterates in  $U_{A_{g_i}}$,
then it needs a fixed number of iterates 
to go from $U_{A_{g_i}}$ to $U_{B_{g_i}}$,
and finally it remains $m_k\,\pi(B_{g_i})$ iterates in
$U_{B_{g_i}}$. In particular, there is a  constant $t\in
\NN$
 independent of $k$ such that
$$
\pi(S_{ik})=m_k\, \pi(B_{g_i})+ \ell_k \, \pi(A_{g_i})+ t.
$$
\end{description}
Moreover, properties (3) and (4) from Proposition~\ref{p.bdf} guarantee that making an arbitrary small perturbation of $g_{ik}$ (preserving properties {\bf S1)} - {\bf S4)}) we can obtain an additional property:
\begin{description}
\item{\bf S5)}
the saddle $S_{ik}$ has the  $V$-$\ut$-biaccumulation property.
\end{description}
\blm\label{l.inequalities}
For every large $k\in \mathbb{N}$ the saddles $S_{ik}$  also have the following properties:
\begin{equation}\label{e.1}
   0<\chi^c(S_{ik})<\frac{1}{2}\chi^c(B_{g_{ik}}),
 \end{equation}
 \begin{equation}\label{e.2}
    \frac{m_k\pi(B_{g_i})}{m_k\pi(B_{g_i})+\ell_k\pi(A_{g_i})+t}>1-C\chi^c(B_{g_{ik}}).
\end{equation}
\elm
\begin{proof}
If $k$ is large enough then $g_{ik}\in \mathcal{U}(g_i).$ Let us
show that for large $k$ the inequality (\ref{e.1}) holds. On the one
hand, from {\bf S2)} and {\bf S4)} we have
$$
0<\chi^c(S_{ik})\le \frac{\log\const}{\pi(S_{ik})}=\frac{\log \const}{m_k\, \pi(B_{g_i})+ \ell_k \,
\pi(A_{g_i})+ t}\to 0 \ \ \ \text{\rm as}\ \ \ \ k\to \infty.
$$
On the other hand, due to (\ref{e.betabis}),  $0<
\dfrac{\log\beta^{\frac{1}{2}}}{\pi(B_{g_{i}})}\le
\chi^c(B_{g_{ik}}).$ This implies that (\ref{e.1}) holds for large
$k$.

Now let us prove that (\ref{e.2}) also holds for large enough $k\in \NN$.
Due to
(\ref{e.lambda}) and
(\ref{e.betabis}), using {\bf S4)}, the central multiplier of
$S_{ik}$ can be estimated from above,
 $$
1<|\lambda^c(S_{ik})|\le |(\beta^2)^{m_k} (\lambda^{\frac{1}{2}})^{\ell_k} T|,
$$
where $T$ is a constant (that does not depend on $k$) which is responsible for the part of the orbit of $S_{ik}$ outside of neighborhoods $U_{A_{g_{i}}}$ and $U_{B_{g_{i}}}$
which does not depend on $k$. 
 This implies that
$$
\ell_k<-\frac{\log T}{\log \lambda^{\frac{1}{2}}}-\frac{m_k\log
\beta^2 }{\log \lambda^{\frac{1}{2}}}.
$$
From this estimate
 for every large $k$ we have
$$
\begin{array}{ll}\label{e.gammaestimate}
    &1-\dfrac{m_k\pi(B_{g_i})}{m_k\pi(B_{g_i})+\ell_k\pi(A_{g_i})+t}=
\dfrac{\ell_k\pi(A_{g_i})+t}{m_k\pi(B_{g_i})+\ell_k\pi(A_{g_i})+t}\le
\dfrac{\ell_k\pi(A_{g_i})+t}{m_k\pi(B_{g_i})}\le\\
 &\qquad \qquad \le
\dfrac{1}{m_k}\left(\dfrac{t}{\pi(B_{g_i})}-\dfrac{\log T}{\log \lambda^{\frac{1}{2}}
}
\cdot\dfrac{\pi(A_{g_{i}})}{\pi(B_{g_i})}
\right)-4\,
\dfrac{\log \beta}{\log \lambda}\cdot\dfrac{\pi(A_{g_{i}})}{\pi(B_{g_i})}=\\
 &\qquad \qquad =\dfrac{1}{m_k}\left(\dfrac{t}{\pi(B_{g_i})}-\dfrac{\log T}{\log
\lambda^{\frac{1}{2}}}
\cdot\dfrac{\pi(A_{g_{i}})}{\pi(B_{g_i})}
\right)-4\, \dfrac{\chi^c(B_{g_{i}})}{\chi^c(A_{g_i})}<8\,\dfrac{\chi^c(B_{g_{i}})}{|\chi^c(A_{g_i})|}.
\end{array}
$$
Finally, due to the choice of the constant $C$ defined
by (\ref{e.C}),
$$
1-\dfrac{m_k\pi(B_{g_i})}{m_k\pi(B_{g_i})+\ell_k\pi(A_{g_i})+t}
<\frac{1}{2}\, C\chi^c(B_{g_{i}})<C\, \chi^c(B_{g_{ik}}).
$$
This completes the proof of the lemma.
\end{proof}

Consider now the set
$$
\mathcal{D'}=\left\{ g_{ik}\  |\  S_{ik} \text{\rm \ \ satisfies
conditions (\ref{e.1}), (\ref{e.2})} \right\} \subset \mathcal{Z}.
$$
By construction, the set $\mathcal{D'}$ is a countable dense subset
of $\mathcal{Z}$. Let us enumerate the elements of
$\mathcal{D'}=\{h_{n_1}\}_{n_1\in \NN}$. Let us also redenote by
$P_{n_1}$ the periodic saddle $S_{ik}$ of the map $h_{n_1}\equiv
g_{ik}$.
 For each $h_{n_1}$ we can apply Proposition~\ref{p.simpleanddense} to the
 heteroclinic $V$-related cycle associated to saddles $P_{n_1}$ and $A_{h_{n_1}}$.
  This gives for each $n_1\in \NN$ an open set $\mathcal{U}_{n_1}\subseteq \mathcal{Z}$
  and a dense countable subset $\widetilde{\mathcal{D}}_{n_1}\subset \mathcal{U}_{n_1}$ such that $h_{n_1}\in \overline{\mathcal{U}}_{n_1}$
   and every $g\in \widetilde{\mathcal{D}}_{n_1}$ has a $V$-related cycle associated with  $A_g$ and a continuation of
   $P_{n_1}$. In order to simplify the notation, we will omit the dependence of the continuation of $P_{n_1}$ on $g$ and
   will write just $P_{n_1}$ instead of ``continuation of $P_{n_1}$".  We can take
    $\mathcal{U}_{n_1}$   small enough to guarantee that for every $g\in \mathcal{U}_{n_1}$ one has
    $0<\chi^c(P_{n_1})<\frac{1}{2}\chi^c(B_g)$. Indeed, due to
    (\ref{e.1}) this inequality holds for $h_{n_1}$. Since
    Lyapunov exponents of a hyperbolic saddle depend continuously on
    a diffeomorphism, the inequality holds also for all $g$ sufficiently $C^1$-close to
    $h_{n_1}$.

Let us now take inductively
$$
\mathcal{Z}_1=\mathcal{U}_{1}, \ \mathcal{Z}_2=\mathcal{U}_{2}\backslash \overline{\mathcal{Z}}_1, \ldots, \mathcal{Z}_{n_1}=\mathcal{U}_{n_1}\backslash \overline{\mathcal{Z}}_{n_1-1}, \ldots
$$
and
$$
\mathcal{D}_{n_1}=\widetilde{\mathcal{D}}_{n_1}\cap
\mathcal{Z}_{n_1}, n_1\in \NN.
$$

We claim that the collection of sets $\{\mathcal{Z}_{n_1}\}_{n_1\in
\NN}$ satisfies the required properties {\bf Z1)} - {\bf Z8)}.
\begin{itemize}
\item
Properties {\bf Z1)} and {\bf Z2)} directly follow from the
construction of $\{\mathcal{Z}_{n_1}\}_{n_1\in \NN}$.
\item
Since the set $\{h_{n_1}\}_{n_1\in \NN}$ is dense in $\mathcal{Z}$,
the union $\cup_{n_1\in \NN}\mathcal{U}_{n_1}$ is dense in
$\mathcal{Z}$, and hence $\cup_{n_1\in \NN}\mathcal{Z}_{n_1}$ is
also dense in $\mathcal{Z}$, so {\bf Z3)} holds.
\item
For each $g\in \mathcal{Z}_{n_1}$ saddle $\{P_{n_1}\}_{n_1\in \NN}$
is $V$-homoclinically related to $B_g$, has real multipliers,
$V$-$\ut$-biaccumulation property, and $\pi(B_g)<\pi(P_{n_1})$, so {\bf
Z4)} holds.
\item
The sets $\mathcal{D}_{n_1}\subset \mathcal{Z}_{n_1}$ were
constructed to satisfy {\bf Z5)}.
\item
Take $g\in \mathcal{Z}_{n_1}$, and denote by $\Gamma$ the part of
the orbit of $P_{n_1}$ that belongs to the neighborhood
$U_{B_g}$. Define the projection
$$
\rho:\Gamma\to
\mathcal{O}(B_g), \quad \rho(x)=\{\text{\rm the point of
$\mathcal{O}(B_g)$ nearest to $x$}\}.
$$
By construction,
$$
\#\Gamma=m_k\pi(B_g)\quad
\mbox{and}
\quad
\#(\mathcal{O}(P_{n_1}))=m_k\pi(B_g)+\ell_k\pi(A_g)+t.
$$
Recall that
 here $k$
and $n_1$ are related due to the enumeration $h_{n_1}=g_{ik}$;
notice that in fact integers $m_k$, $\ell_k$, and $t$ depend also on
the index $i$, but our notations do not reflect this dependence. Now
{\bf Z6)} follows from the inequality (\ref{e.2}) and the choice
of $\gamma$ in (\ref{e.gamma}). 
\item
The value of $\gamma$ could be taken arbitrary small; in particular
{\bf Z7)} can be satisfied by the choice of sufficiently small
$\gamma$ (that choice was explicitly specified in \eqref{e.gamma}).
\item
The last property {\bf Z8)} follows directly from the inequality
(\ref{e.1}).
\end{itemize}
Finally, assume that the sets $\{\mathcal{Z}_{n_1, \ldots ,n_N},
n_i\in \mathbb{N}\}$ were constructed. Take one of these sets, say,
$\mathcal{Z}_{n_1, \ldots ,n_N}$. Exactly the same arguments that we
used to construct the sets $\mathcal{Z}_{n_1}\subset \mathcal{Z},
n_1\in \mathbb{N}$, can be now used to construct the sets
$\mathcal{Z}_{n_1, \ldots, n_N, n_{N+1}}\subseteq \mathcal{Z}_{n_1,
\ldots, n_N}, n_{N+1}\in \mathbb{N}$. By induction,
Proposition~\ref{p.conditions}
 follows.
\end{proof}

\vskip 1cm

\flushleft{\bf Lorenzo J. D\'\i az}
 \ \ (lodiaz@mat.puc-rio.br)\\
Departamento de  Matem\'{a}tica, PUC-Rio \\ Marqu\^{e}s de S. Vicente 225\\
22453-900 Rio de Janeiro RJ \\ Brazil
\medskip

\flushleft
{\bf Anton Gorodetski}  \ \  (asgor@math.uci.edu)\\
 Department of Mathematics\\
University of California, Irvine\\
 Irvine, CA 92697\\
USA 

\end{document}